\RequirePackage{fix-cm}
\documentclass[smallextended]{svjour3}       % onecolumn (second format)
\smartqed  % flush right qed marks, e.g. at end of proof
\usepackage{amsmath,amssymb,amsfonts,graphicx,float,multicol,fleqn}
\usepackage{multirow}
\usepackage{graphicx}
\usepackage{array}
\usepackage{booktabs}
\usepackage{subfigure}

\begin{document}

\title{A matrix formulation of the Tau method for the numerical solution of non-linear problems}

\titlerunning{Tau method for the numerical solution of non-linear problems}        % if too long for running head

\author{ Kourosh Parand \and Amin Ghaderi  \and Mehdi Delkhosh \and Reza Pourgholi}


\institute{K. Parand \at
              Department of Computer Sciences, Shahid Beheshti University, G.C., Tehran, Iran. \\
              Department of Cognitive Modelling, Institute for Cognitive and Brain Sciences, Shahid Beheshti University, G.C, Tehran, Iran. \\
              \email{k\_parand@sbu.ac.ir}           %  \\
           \and
           A. Ghaderi\at
            Department of Cognitive Modelling, Institute for Cognitive and Brain Sciences, Shahid Beheshti University, G.C, Tehran, Iran.\\
              \email{amin.g.ghaderi@gmail.com}
	\and
 	M. Delkhosh\at
	Department of Mathematics and Computer Science, Islamic Azad University, Bardaskan Branch, Bardaskan, Iran.\\
  \email{ mehdidelkhosh@yahoo.com}
           \and
           R. Pourgholi \at
            School of  Mathematics and Computer  Sciences, Damghan University, Damghan, Iran. \\
              \email{pourgholi@du.ac.ir}
}

\date{Received: date / Accepted: date}

\maketitle

\begin{abstract}
The purpose of this research is to propose  a new approach named the shifted Bessel Tau (SBT) method for solving higher-order ordinary differential equations (ODE). The operational matrices of derivative, integral and product of shifted Bessel polynomials on the interval $[a,b]$ are calculated. These matrices together with the Tau method are utilized to reduce the solution of the higher-order ODE to the solution of a system of algebraic equations with unknown Bessel coefficients. The comparisons between the results of the present work and other the numerical method are shown that the present work is computationally simple and highly accurate.

\end{abstract}

\section{Introduction}\label{sec1}
Differential equations and their solutions play a major role in science and engineering. A physical event can be modeled by the differential equation, an integral equation or an integro-differential equation or a system of these equations. Since few of these equations cannot be solved explicitly, it is often necessary to resort to the numerical techniques which are appropriate combinations of numerical integration and interpolation. The solution of this equation occurring in physics, biology and engineering are based on numerical methods. One of the powerful method for solving the ordinary differential equations, integral and integro-differential equations are spectral methods in both terms of accuracy and simplicity. Spectral methods, in the context of numerical schemes for differential equations, generically belong to the family of weighted residual methods (WRMs) \cite{refIntroduction01}. WRMs represent a particular group of approximation techniques, in which the residuals (or errors) are minimized in a certain way and thereby leading to specific methods, including Galerkin, Petrov-Galerkin, collocation and Tau formulations \cite{refIntroduction02,refIntroduction03,refIntroduction04,refIntroduction05}. The basis of spectral method is always polynomials (orthogonal polynomials). These types of basis for approximation of enough smooth functions in the finite domain leads to the exponential convergence.\\\
In Tau method, the operational matrices of derivative, integral and product of basis functions are obtained  and then these matrices together to reduce the solution of these problems to the solution of a system of algebraic equations. Tau method is based on expanding the required approximate solution as the elements of a complete set of basis functions. Polynomial series and orthogonal functions have received considerable attention in dealing with various problems of differential, integral and integro-differential equations. Polynomials are incredibly useful mathematical tools as they are simply defined and can be calculated quickly on computer systems and represent a tremendous variety of functions. They can be differentiated and integrated easily, and can be pieced together to form spline curves that can approximate any function to any accuracy desired. Also, the main characteristic of this technique is that it reduces these problems to those of solving a system of algebraic equations, thus greatly simplifies the problems.\\
In recent years, the differential, integral and integro-differential equations have been solved using the homotopy perturbation method \cite{refIntroduction12,refIntroduction13}, the fractional order of rational Bessel functions collocation (FRBC) method \cite{Parand.Ghaderi.Electonical}, the radial basis functions method \cite{refIntroduction14}, the collocation method \cite{refIntroduction15}, the Homotopy analysis method \cite{refIntroduction16,refIntroduction17}, the Tau method \cite{refIntroduction18,refIntroduction19}, the Variational iteration method  \cite{refIntroduction20,refIntroduction21}, the Legendre matrix method \cite{refIntroduction10,refIntroduction22}, the Haar Wavelets operational matrix \cite{refIntroduction23}, the Shifted Chebyshev direct method \cite{refIntroduction24}, the generalized fractional order of the Chebyshev collocation approach \cite{refIntroduction42}, the Legendre Wavelets operational matrix \cite{refIntroduction25}, Sine-Cosine Wavelets operational matrix \cite{refIntroduction26}, the operational matrices of Bernstein polynomials \cite{refIntroduction27}. Razzaghi et al. in \cite{refIntroduction28,refIntroduction29,refIntroduction30} presented the integral and product operational matrices based on the Fourier, the Taylor and the shifted-Jacobi series. Parand et al. in \cite{refIntroduction30_1,refIntroduction30_2,Parand.Ghaderi.European,Parand.Ghaderi.SeMA} applied the spectral method to solve nonlinear ordinary differential equations on semi-infinite intervals. Their approach was based on a rational Tau method. They obtained the operational matrices of derivative and product of rational Chebyshev and Legendre functions and then applied these matrices together with the Tau method to reduce the solution of these problems to the solution of a system of algebraic equations. Doha \cite{refIntroduction31} has derived the shifted Jacobi operational matrix of fractional derivatives which is applied together with the spectral Tau method for the numerical solution of dynamical systems.Yousefi et al. in \cite{refIntroduction32,refIntroduction33,refIntroduction34,refIntroduction35} have presented Legendre wavelets and Bernstein operational matrices and used them to solve miscellaneous systems. Also, recently, in \cite{refIntroduction36} a new shifted-Chebyshev operational matrix of fractional integration of arbitrary order is introduced and applied together with the spectral Tau method for solving linear fractional differential equations (FDEs). In \cite{refIntroduction37,refIntroduction38} the operational matrix form in a Hybrid of block-pulse functions and another set of functions like Taylor series, Legendre and Chebyshev has been used for finding the solution of the various classes of dynamical systems. Lakestani \cite{refIntroduction39} constructed the operational matrix of fractional derivative of order $\alpha$ in the Caputo sense using the linear B-spline functions. In \cite{refIntroduction40}, a general formulation for the d dimensional orthogonal functions and their derivative and product matrices has been presented. Ordinary operational matrices have been utilized together with the Tau method to reduce the solution of partial differential equations (PDEs) to the solution of a system of algebraic equations. Recently in \cite{refIntroduction41}, a class of two dimensional nonlinear Volterra's integral equations has been solved using operational matrices of Legendre polynomials. The operational matrices of integration and product together with the collocation points have been utilized to reduce the solution of the integral equation to the solution of a nonlinear algebraic equation system.

The remainder of the paper is organized as follows: In section \ref{sec2}, the function approximation and convergence of the method by using Bessel polynomials  are illustrated . In sections  \ref{sec3}, the general procedure for formulation of operational matrices of integration, differentiation, and product are explained  and are proved. In section \ref{sec4}, our numerical findings and demonstrate the validity, accuracy and applicability of the operational matrices by considering numerical examples  are reported. Also a conclusion is given in the last section. 

\section{Bessel Polynomials}\label{sec2} The Bessel functions arise in many problems in physics possessing cylindrical symmetry, such as the vibrations of circular drumheads and the radial modes in optical fibers. Bessel functions are usually defined as a particular solution of a linear differential equation of the second order which known as Bessel's equation. 

Bessel differential equation of order $n\in\mathbb{R}$ is:
\begin{eqnarray}\label{equation01}
x^2\frac{d^2y(x)}{dx^2}+x\frac{dy(x)}{dx}+(x^2-n^2) y(x)=0,~~x\in(-\infty,\infty).
\end{eqnarray}

One of the solutions of Eq. (\ref{equation01}) by applying the method of Frobenius is as follows \cite{refBesselPol03}:
\begin{eqnarray}\label{equation02}
J_{n}(x)=\sum_{r=0}^{ \infty}\frac{(-1)^r}{r!(n+r)!}(\frac{x}{2})^{2r+n},
\end{eqnarray}
where series (\ref{equation02})  is convergent for all $x\in(-\infty,\infty)$.

Bessel functions and polynomials are used to solve the more number of problems in physics, engineering, mathematics, and etc., such as Blasius equation, Lane-Emden equations, integro-differential equations of the fractional order, unsteady gas equation, systems of linear Volterra integral equations, high-order linear complex differential equations in circular domains, systems of high-order linear Fredholm integro-differential equations, nolinear Thomas-Fermi on semi-infinite domain, etc. \cite{refBesselPol04,refBesselPol05,refBesselPol09,refBesselPol10,refBesselPol11,refBesselPol12,refBesselPol13,refBesselPol14,refBesselPol15,refBesselPol16,Parand.Ghaderi.Electonical}.

Bessel polynomials has been introduced as follows \cite{refBesselPol12}:
\begin{eqnarray}
\nonumber B_{n}(x)=\sum_{r=0}^{[\frac{N-n}{2}]}\frac{(-1)^r}{r!(n+r)!}(\frac{x}{2})^{2r+n},~~x\in[0,1].
\end{eqnarray}
where $n\in\mathbb{N}$, and $N$ is the number of basis of Bessel polynomials.\\
Now, We introduce new shifted Bessel equation in the interval $[a,b]$ by:
\begin{eqnarray}
\nonumber Q_{n}(x)=B_n(\frac{x-a}{b-a}),~~x\in[a,b],~~n=0, 1, 2, ..., N.
\end{eqnarray}
where $a, b \in \mathbb{R}$.

\subsection{Approximation of functions}
Let us define $\Gamma=\{x|~a\leq x \leq b \}$ and \\
$L^{2}_{w}(\Gamma)=\{~v :\Gamma \rightarrow \mathbb{R}| v$ is measurable and $\parallel v \parallel_{w} < \infty \}$, where
$$\parallel v \parallel_{w}=\left(\int^{b}_{a}|v(x)|^{2}w(x)dx\right)^{1/2},$$
with $w(x)=\frac{1}{b-a}$, is the norm induced by inner product of the space $L^{2}_{w}(\Gamma)$ as follows:
$$\langle v,u\rangle_{w}=\int^{b}_{a}{v(x)u(x)w(x)}dx.$$

Now, suppose that\\
$$\mathfrak{Q}_N=~span\{Q_{0}(x), Q_{1}(x),\dots, Q_{N}(x)\},$$
where $\mathfrak{Q}_N$ is a finite-dimensional subspace of $L^{2}_{w}(\Gamma)$, dim ($\mathfrak{Q}_N) = N+1$, so $\mathfrak{Q}_N$ is a closed subspace of $L_{w}^{2}(\Gamma)$. Therefore, $\mathfrak{Q}_N$ is a complete subspace of $L_{w}^{2}(\Gamma)$. Assume that $f(x)$ be an arbitrary element in $L_{w}^{2}(\Gamma)$. Thus $f(x)$ has a unique best approximation in $\mathfrak{Q}_N$ subspace, say $\hat{b}\in \mathfrak{Q}_N$, that is
\begin{eqnarray}
\nonumber\exists~ \hat{b(x)}\in\mathfrak{Q}_N, ~~~ \forall ~b(x)\in \mathfrak{Q}_N,~~\parallel f(x)-\hat{b(x)}\parallel_{w} \leq \parallel f(x)-b(x)\parallel_{w}.
\end{eqnarray}
Notice that we can write $b(x)$ function as a linear combination of the basis vectors of $\mathfrak{Q}_N$ subspace. 

We know function of $f(x)$  can be expanded by $N+1$ terms of Bessel polynomials as: 
\begin{eqnarray}
f(x)=f_{N}(x)+R(x),
\nonumber\end{eqnarray}
that is
\begin{equation}\label{equation03}
f_{N}=\sum^{N}_{n=0}{a_{n}Q_{n}(x)}=A^{T}Q(x),
\end{equation}
where $Q(x)=[Q_{0}(x), Q_{1}(x),\cdots, Q_{N}(x)]^{T}$ and $R(x)\in\mathfrak{Q}^{\perp}_N$ that ${\mathfrak Q}^{\perp}_N$ is the orthogonal complement. So  $f(x)-f_{N}(x)\in\mathfrak{Q}_N^{\perp}$ and $b(x)\in\mathfrak{Q}_N$ are orthogonal which we denote it by:
\begin{eqnarray}
\nonumber \big(f(x)-f_{N}(x)\big)\perp b(x),
\end{eqnarray}
thus $f(x)-f_{N}(x)$ vector is orthogonal over all of basis vectors of $\mathfrak{Q}_N$ subspace as:
\begin{eqnarray}
\nonumber \langle f(x)-f_{N}(x),Q_{i}(x)\rangle_{w}=\langle f(x)-A^{T}Q(x),Q_{i}(x)\rangle_{w}=0,~i=0, 1,\cdots, N,
\end{eqnarray}
hence
\begin{eqnarray}
&&\nonumber\langle f(x)-A^{T}Q(x),Q^{T}(x)\rangle_{w}=0,
\end{eqnarray}
therefore A can be obtained by
\begin{eqnarray}
&&\nonumber\langle f(x),Q^{T}(x)\rangle_{w}=\langle A^{T}Q(x),Q^{T}(x)\rangle_{w},\\
\nonumber\\
\nonumber &&A^{T}=\langle f(x),Q^{T}(x)\rangle_{w}\langle Q(x),Q^{T}(x)\rangle_{w}^{-1},
\end{eqnarray}
where
\begin{eqnarray}\label{equation04}
&& K=\langle Q(x),Q^{T}(x)\rangle_{w}=\int_{a}^{b}{Q(x)Q(x)^{T}}w(x)dx,
\end{eqnarray}
is an $(N +1)×(N +1)$ matrix and is said dual matrix of $Q(x)$.
\subsection{Convergence of the method}
The following theorem shows that by increasing $N$, the approximation solution $f_N(x)$  is convergent to $f(x)$ exponentially.\\\\
\textbf{Theorem 1.} (\textit{Taylor's formula}) Suppose that $f(x) \in C[0,b]$ and $D^{k}f(x) \in C[0,b]$, where $k=0,1,...,N$. Then we have 
\begin{equation}\label{eqn3}
f(x) = \sum_{i=0}^{N} \frac{x^{i}}{i!}D^{i}f(0^+)+\frac{x^{N+1}}{(N+1)!}D^{N+1}f(\xi),
\end{equation}
with $0\leqslant\xi \leqslant x,~ \forall x\in [0,b]$. And thus
\begin{equation}\label{eqn4}
| f(x) - \sum_{i=0}^{N} \frac{x^{i}}{i!}D^{i}f(0^+) | \leqslant M \frac{x^{N+1}}{(N+1)!},
\end{equation}
where $M \geqslant |D^{N+1}f(\xi)|$.\\
\textbf{Proof:} See Ref. \cite{Odibat.Momani}.

\textbf{Theorem 2} 
Suppose that $D^{k}f(x)\in C[0,b]$ for $k=0,1,...,N,$ and $\mathfrak{Q}_{N}$ is the subspace generated by $\{Q_0(x),~Q_1(x),...,~Q_{N}(x)\}$. If $f_N(x)=F^TQ(x)$ is the best approximation to $f$ from $\mathfrak{Q}_{N}$ , then the error bound is presented as follows
\begin{eqnarray}\label{equation05}
\parallel f(x)-f_N(x)\parallel_w \leqslant \frac{M}{(N+1)!}\big (\frac{b^{2N+2}}{(2N+3)}\big)^{\frac{1}{2}},
\end{eqnarray}where $M \geqslant|D^{N+1}f(x)|, ~~x\in [0,b].$\\
\textbf{Proof.} By theorem 1, we have $y(x)= \sum_{i=0}^{N} \frac{x^{i}}{i!}D^{i}f(0^+)$ and 
$$| f(x)-y(x) | \leqslant M \frac{x^{N+1}}{N+1!}.$$
Since the best approximation to $f(x)$ from $\mathfrak{Q}_N$ is $F^TQ(x)$, and $y(x)\in  \mathfrak{Q}_N$, thus
\begin{eqnarray}
\parallel f(x)-f_N(x)\parallel_w^2 \leqslant \parallel f(x)-y(x)\parallel_w^2 \leqslant \frac{M^2}{(N+1)!^2}\int_0^bx^{2N+2}~w(x)~dx, \nonumber
\end{eqnarray}
where the weight function is $w(x)=\frac{1}{b}$, thus by integration of above equation, Eq. (\ref{equation05}) can be proved $\blacksquare$.

\section{Operational matrices}\label{sec3}
Now, we can show $Q(x)$ as metrix product of $\textbf{Y}$ and $X_{s}(x)$ as follows:
\begin{eqnarray}\label{equation06}
Q(x)=\textbf{Y}X_{s}(x),
\end{eqnarray}
where $Q(x)=[Q_0(x),Q_1(x),...,Q_{N}(x)]^T$ and $X_s(x)=[1,\big(\frac{x-a}{b-a}\big),\big(\frac{x-a}{b-a}\big)^2,...,\big(\frac{x-a}{x-b}\big)^N]^T$.

There are two different forms for $\textbf{Y}$. If $N$ is odd \cite{refBesselPol10,refBesselPol11,refBesselPol12,refBesselPol13,refBesselPol14,refBesselPol15}:

\textbf{Y}=$ \left[ \begin{array}{cccccc}  
\frac{1}{0!0!2^0} &0                             & \frac{-1}{1!1!2^2}& \cdots & \frac{(-1)^{\frac{N-1}{2}}}{(\frac{N-1}{2})!(\frac{N-1}{2})!2^{N-1}} &0\\  
0                             & \frac{1}{0!1!2^1}&0                              & \cdots & 0                                                                                                                  & \frac{(-1)^{\frac{N-1}{2}}}{(\frac{N-1}{2})!(\frac{N+1}{2})!2^{N}} \\
0                             &0                             & \frac{1}{0!2!2^2} & \cdots &  \frac{(-1)^{\frac{N-3}{2}}}{(\frac{N-3}{2})!(\frac{N+1}{2})!2^{N-1}}&0\\ 
\vdots                     &\vdots                    &\vdots                     & \cdots & \vdots                                                                                                           & \vdots\\ \\
0                             &0                            &0                              & \cdots &\frac{1}{ 0!(N-1)!2^{N-1}}                                                                           &0 \\ 
0                             &0                            &0                              & \cdots & 0                                                                                                                    & \frac{1}{0!N!2^{N}}
\end{array} \right]$

And if $N$ is even:

\textbf{Y}=$ \left[ \begin{array}{cccccc}  
\frac{1}{0!0!2^0} &0                             & \frac{-1}{1!1!2^2}& \cdots &0											   & \frac{(-1)^{\frac{N}{2}}}{(\frac{N}{2})!(\frac{N}{2})!2^{N}} \\  
0                             & \frac{1}{0!1!2^1}&0                              & \cdots & \frac{(-1)^{\frac{N-2}{2}}}{(\frac{N-2}{2})!(\frac{N}{2})!2^{N-1}} & 0                                                                                                   \\
0                             &0                             & \frac{1}{0!2!2^2} & \cdots & 0 										    & \frac{(-1)^{\frac{N-2}{2}}}{(\frac{N-2}{2})!(\frac{N+2}{2})!2^{N}}\\ 
\vdots                     &\vdots                    &\vdots                     & \cdots & \vdots 										    & \vdots \\ 
0                             &0                            &0                              & \cdots &\frac{1}{ 0!(N-1)!2^{N-1}} 							    &0	 \\ 
0                             &0                            &0                              & \cdots & 0										    &\frac{1}{0!N!2^{N}}\\
\end{array} \right]$

In this step, $X_s(x)$ as metrix product of $\textbf{S}$ and $X(x)$  are expaned as follows:
\begin{eqnarray}\label{equation06}
X_s(x)=\textbf{S}X(x),
\end{eqnarray}
which $X(x)=[1,x,...,x^N]^T$ and 
\begin{eqnarray*}
\textbf{S}=\begin{bmatrix}
1&0&0&\cdots&0\cr
-\frac{a}{b-a}&\frac{1}{b-a}&0&\cdots&0\cr
\frac{a^2}{(b-a)^2}&-\frac{2a}{(b-a)^2}&\frac{1}{(b-a)^2}&\cdots&0\cr
\vdots&\vdots&\vdots&\vdots&\vdots\cr
\frac{\dbinom{N}{0}}{\big(b-a\big)^N\big(-a\big)^{-N}}&\frac{\dbinom{N}{1}}{\big(b-a\big)^N\big(-a\big)^{-N+1}}&\frac{\dbinom{N}{2}}{\big(b-a\big)^N\big(-a\big)^{-N+2}}&\cdots&\frac{1}{\big(b-a\big)^N}\cr
\end{bmatrix}
\end{eqnarray*}
Therfore, $Q(x)$ will be rewritten as follows:
\begin{eqnarray}\label{equation06}
Q(x)=\textbf{M}X(x),
\end{eqnarray}
where $\textbf{M}=\textbf{Y} \textbf{S}$.\\
We know a triangular matrix is invertible if and only if all its diagonal entries are nonzero, and  the product of two invertible matrices also invertible, so
\begin{eqnarray}\label{equation07}
X(x)=\textbf{M}^{-1}Q(x).
\end{eqnarray}
This shows that we can represent  $X(x)$ as metrix product of $\textbf{M}^{-1}$ and $Q(x)$.
\subsection{The operational matrix of derivative}
Let start this section by introducing the operational matrix of derivative of $X(x)$ vector as follows:
\begin{eqnarray}\label{equation08}
\frac{d}{dx}X(x)=P~X(x),
\end{eqnarray}
where
\[ P_{i,j}= \left\{
\begin{array}{ll}
      i,& i-j=1\\
      0,& otherwise\\
\end{array} ~~for~i,j=0, 1,...,N~,
\right. \]
The matrix $\textbf{D}$ is named as the (N + 1)$\times$(N + 1) operational matrix of derivative if and only if
\begin{eqnarray} \label{equation09}
\frac{d}{dx}Q(x)=\textbf{D}~Q(x).
\end{eqnarray}
\textbf{Lemma 1.} Suppose that   $\textbf{M}$, $\textbf{M}^{-1}$ and $\textbf{D}$ are (N + 1)$\times$(N + 1) matrices which satisfy respectively  Eqs. (\ref{equation06}), (\ref{equation07}) and (\ref{equation09}). Then
\begin{eqnarray}\label{equation10}
\textbf{D}=\textbf{M}P\textbf{M}^{-1}.
\end{eqnarray}
\textbf{Proof:} By using expansion of $\frac{d}{dx}Q(x)$, we have
\begin{eqnarray}
\nonumber&&\frac{d}{dx}Q(x)=\frac{d}{dx}\textbf{M}X(x)=\textbf{M}\frac{d}{dx}X(x)=\textbf{M}PX(x)=\textbf{M}P\textbf{M}^{-1}Q(x)\\
\nonumber &&\Rightarrow \textbf{D}=\textbf{M}P\textbf{M}^{-1}~  \blacksquare.
\end{eqnarray}
\textbf{Remark:}
If $\textbf{D}$ be a first order operational matrix of drivative of $Q(x)$ then $k$th order operational matrix of drivative of $Q(x)$ is as follows:
\begin{eqnarray}\label{equation11}
\frac{d^k}{dx^k}Q(x)=\textbf{D}^kQ(x).
\end{eqnarray}
\subsection{Operational matrix of integral}
The integral operational matrix for basis vector $X(x)$ operates as:
\begin{eqnarray}\label{equation12}
\int_{a}^{x}X(t)dt\simeq LX(x),
\end{eqnarray}
where
\[ L_{i,j}= \left\{
\begin{array}{lll}
      \frac{1}{i+1}&j-i=1~and~i\neq N\\
      -\frac{a^{i+1}}{i+1}&j=0~and~i\neq N\\
     r_{j}&i=N\\
      0&otherwise\\
\end{array}~~for~i,j=0, 1,...,N~,
\right. \]
where $\nonumber R=[r_0,r_1,...,r_{m}]^T$.

 We just need to approximate $\int_{a}^{x}t^{N}dt=\frac{x^{N+1}}{N+1}-\frac{a^{N+1}}{N+1}$. 

According to the relation (\ref{equation03}), we have
\begin{eqnarray}\label{equation13}
\frac{x^{N+1}}{N+1}-\frac{a^{N+1}}{N+1}\simeq \sum_{n=0}^{N} a_n ~Q_n(x)=R^TQ(x),
\end{eqnarray}
thus
\begin{eqnarray}
R^T \simeq  \langle \frac{x^{N+1}}{N+1},Q^T(x) \rangle _w ~ \langle Q(x),Q^T(x) \rangle _w^{-1},~~for~w(x)=\frac{1}{b-a}. \nonumber
\end{eqnarray}

The matrix $I$ is named as the (N + 1)$\times$(N + 1) operational matrix of the derivative of the vector Q(x) if and only if
\begin{eqnarray} \label{equation14}
\int_{0}^{x}Q(t)dt\simeq \textbf{I}Q(x).
\end{eqnarray}
\textbf{Lemma 2.} Suppose that   $\textbf{M}$, $\textbf{M}^{-1}$ and $L$ are  the (N + 1)$\times$(N + 1) matrices which satisfy respectively  Eqs (\ref{equation06}), (\ref{equation07}) and (\ref{equation14}), then
\begin{eqnarray}\label{equation15}
\textbf{I}=\textbf{M}L\textbf{M}^{-1}.
\end{eqnarray} 
\textbf{Proof:} By using expansion of $\int_{a}^{x}Q(t)dt,$ we have
\begin{eqnarray}
\nonumber&&\int_{0}^{x}Q(t)dt=\int_{0}^{x}\textbf{M}X(t)dt=\textbf{M}\int_{0}^{x}X(t)dt\simeq \textbf{M}LX(x)=\textbf{M}L\textbf{M}^{-1}Q(x) \\
\nonumber &&\Rightarrow \textbf{I}=\textbf{M}L\textbf{M}^{-1}~ \blacksquare.
\end{eqnarray}

\subsection{The operational matrix of product}
The following property of the product of two $X(x)$ vectors will also be applied as follows:
\begin{eqnarray} \label{equation16}
X(x)X(x)^TV\simeq \widetilde{V}X(x)^T,
\end{eqnarray}
where $V = [v_{0}, v_{1}, ..., v_{N}]$. According to the relation (\ref{equation03}), we have
\begin{eqnarray}
\widetilde{V} \simeq  \langle X(x),X^T(x) \rangle _w ~ \langle X(x),X^T(x) \rangle _w^{-1},~~for~w=\frac{1}{b-a}. \nonumber
\end{eqnarray}
The matrix $\widetilde{\textbf{C}}$ is named as (N + 1)$\times$(N + 1) operational matrix of product of Q(x) if and only if
\begin{eqnarray} 
\nonumber Q(x)Q(x)^TC\simeq \widetilde{\textbf{C}}Q(x)^T,
\end{eqnarray}
where $C = [c_{0}, c_{1}, ..., c_{N}]$.

\textbf{Lemma 3.} Suppose that   $\textbf{M}$, $\textbf{M}^{-1}$ and $\widetilde{V}$ are (N + 1)$\times$(N + 1) matrices which satisfy respectively  Eqs. (\ref{equation06}), (\ref{equation07}) and (\ref{equation16}). Then
\begin{eqnarray}\label{equation17}
\widetilde{\textbf{C}}=\textbf{M}\widetilde{V}\textbf{M}^{-1}.
\end{eqnarray} 
\textbf{Proof:} Assume that $V=\textbf{M}^{T}C$ so by using expansion of $ Q(x)Q(x)^TC,$ we have
\begin{eqnarray}
\nonumber&& Q(x)Q(x)^TC=\textbf{M}X(x)\textbf{M}X(x)C=\textbf{M}X(x)X(x)^{T}\textbf{M}^{T}C\\
\nonumber&&~~~~~~~~~~~~~~~~~=\textbf{M}X(x)X(x)V\simeq \textbf{M}\widetilde{V}X(x)=\textbf{M}\widetilde{V}\textbf{M}^{-1}Q(x) \\
\nonumber &&\Rightarrow \widetilde{\textbf{C}}=\textbf{M}\widetilde{V}\textbf{M}^{-1}~ \blacksquare.
\end{eqnarray}

\section{Illustrative examples}\label{sec4}
In the present paper, to illustrate the effectiveness of Tau method, several test examples are carried out. A comparison of our results with those obtained by other methods reveals that our methods are very effective and convenient.
Now we will demonstrate the use of the SBT method in several engineering problems.\\\\
\textbf{Example 1: }
We consider the governing equation for unsteady two-dimensional flow and heat transfer of a viscous fluid \cite{refExapmle01} which the equation is convertible to the following system of the nonlinear equations:
\begin{eqnarray}
&&f'''' -S\left(x f''' +3f'' -2f f''  \right) -G^2f''=0,\label{equation18}\\
&&\theta'' +S~Pr\left(2 f \theta' -x \theta'  \right) +Pr~Ec (f''^2  +12\delta^2 f'^2  )=0, \label{equation19}
\end{eqnarray}
with the boundary conditions
\begin{eqnarray}\label{equation20}
&&f(0)=A,~~~f'(0)=0,~~~\theta(0)=1,\nonumber \\
&&f(1)=\frac{1}{2},~~~f'(1)=0,~~~\theta(1)=0.
\end{eqnarray}
To solve this example, we approximate $f(x)$ and $\theta(x)$  functions by $N+1$ basis of the shifted Bessel polynomials on the interval $[0,1]$ as bellow:
\begin{eqnarray}
&&f(x)\simeq \sum_{i=0}^{N} {a_{i}Q_{i}(x)}=F^{T}Q(x),\label{equation21}\\
&&\theta(x)\simeq \sum_{i=0}^{N} {c_{i}Q_{i}(x)}=C^{T}Q(x),\label{equation22}
\end{eqnarray}
where $F=[a_0,a_1,...,a_{N}]^T$ and $C=[c_0,c_1,...,c_{N}]^T$.

\begin{figure}[!ht]
\centering
\subfigure[Effect of A on $f'(x)$]{
\includegraphics*[width=5cm]{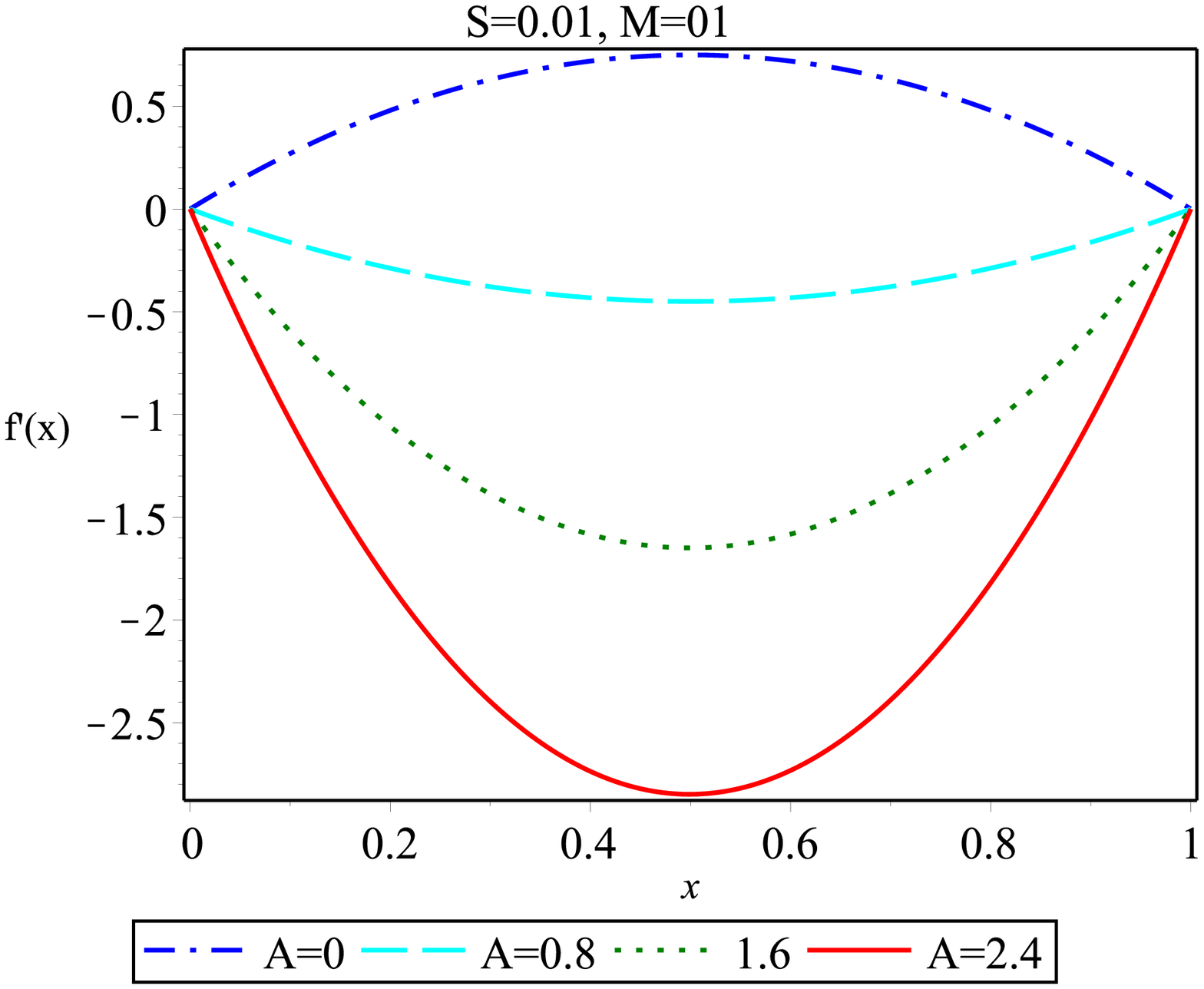}}
\hspace{3mm}
\subfigure[Effect of S on $f'(x)$]{
\includegraphics*[width=5cm]{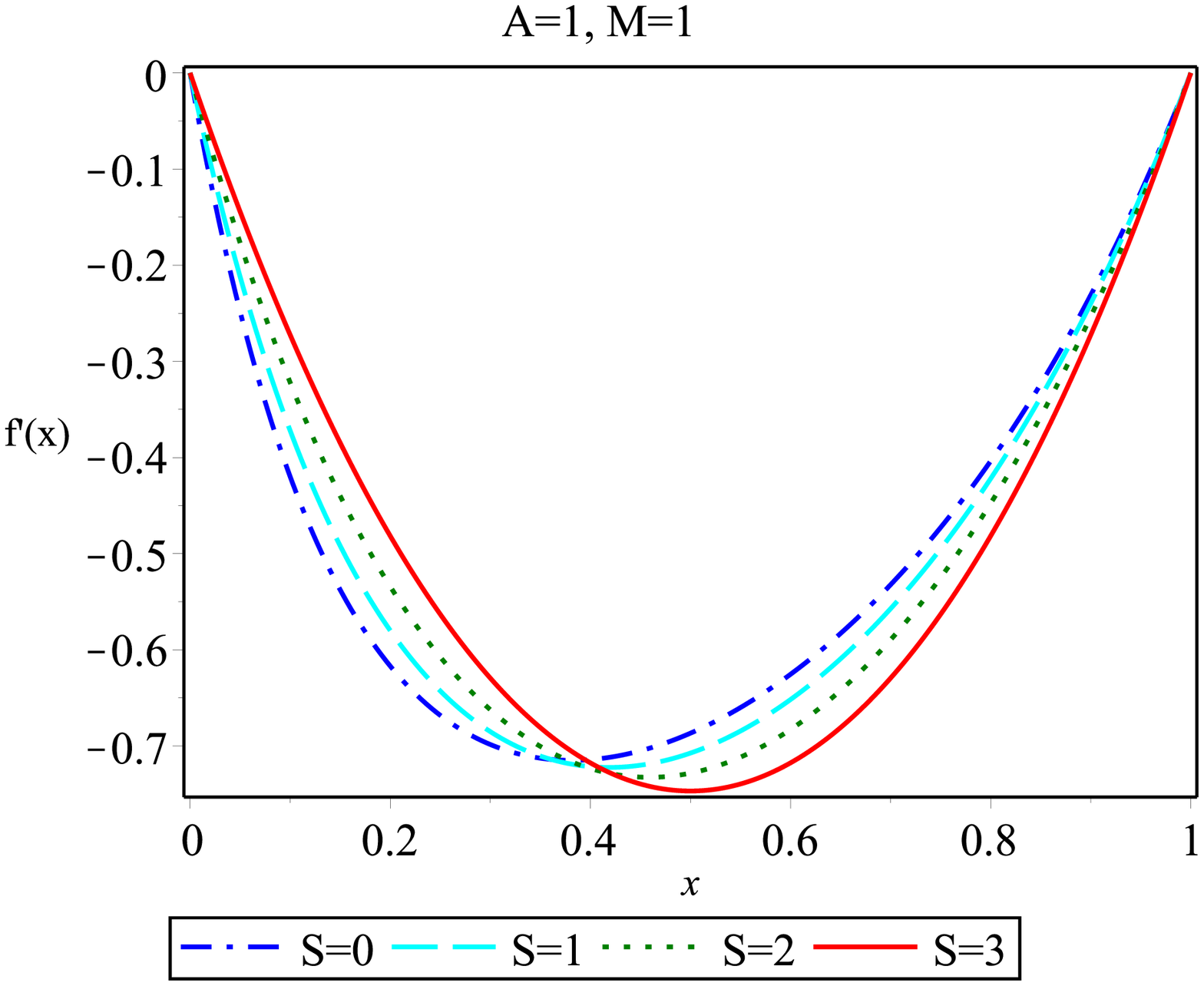}}
\caption{Effect of A and S on $f'(x)$ for the suction cases, Example 1.}
\label{FSection}
\end{figure}

\begin{figure}[!ht]
\centering
\subfigure[Effect of A on $f'(x)$]{
\includegraphics*[width=5cm]{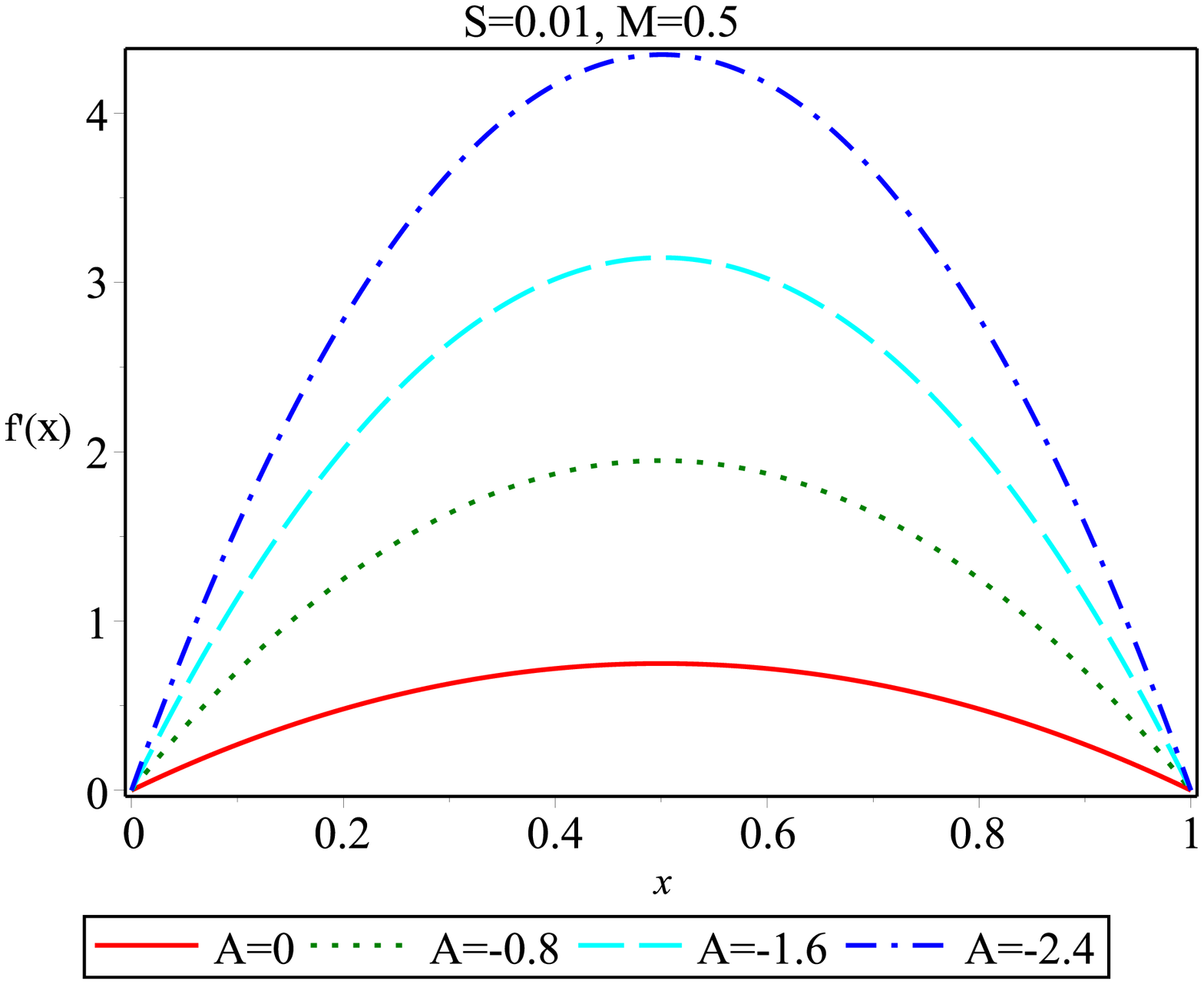}}
\hspace{3mm}
\subfigure[Effect of S on $f'(x)$]{
\includegraphics*[width=5cm]{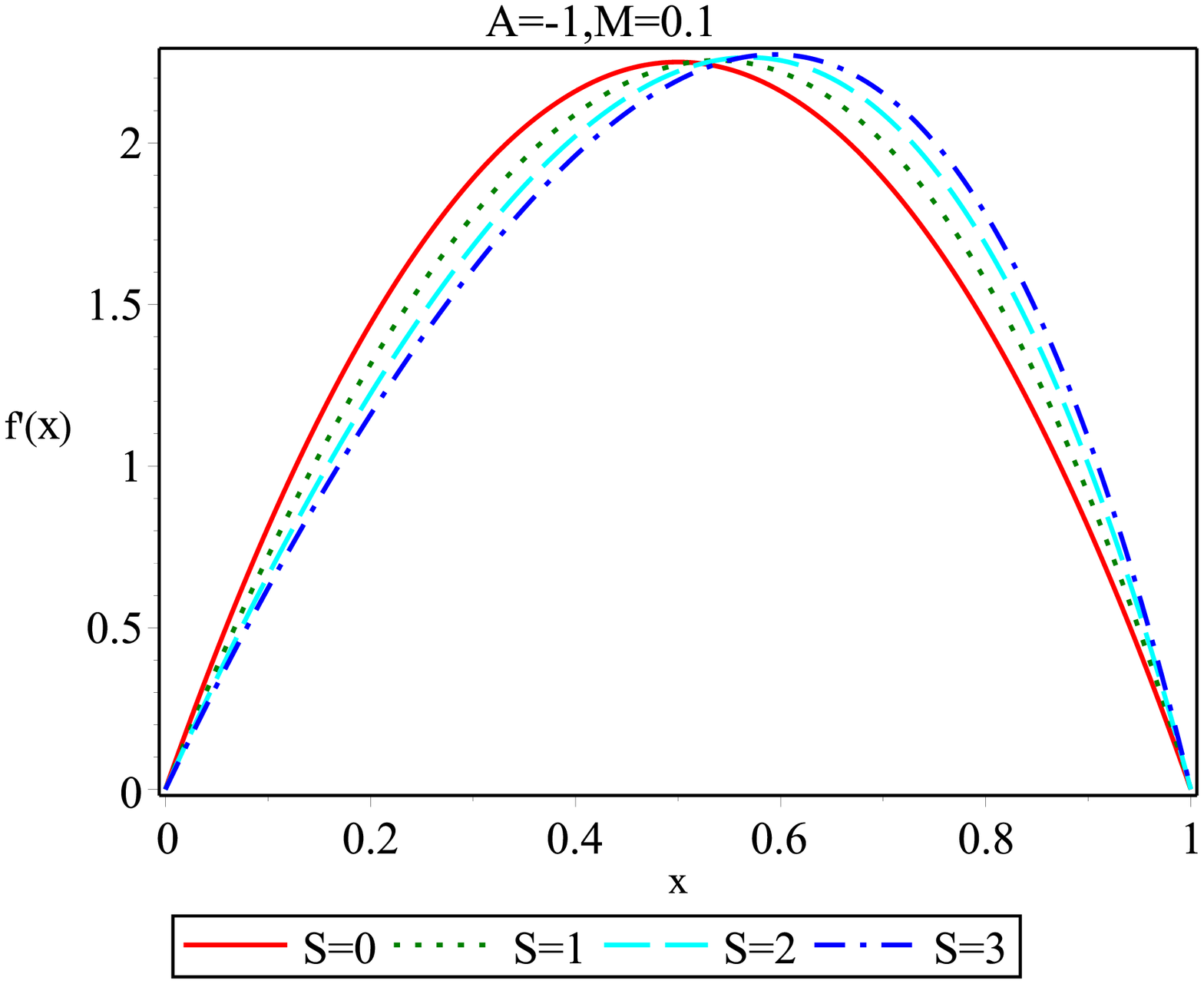}}
\caption{Effect of A and S on $f'(x)$ for the injection cases, Example 1.}
\label{FInjection}
\end{figure}

\begin{figure}[!ht]
\centering
\subfigure[Effect of Pr on $\theta(x)$]{
\includegraphics*[width=5cm]{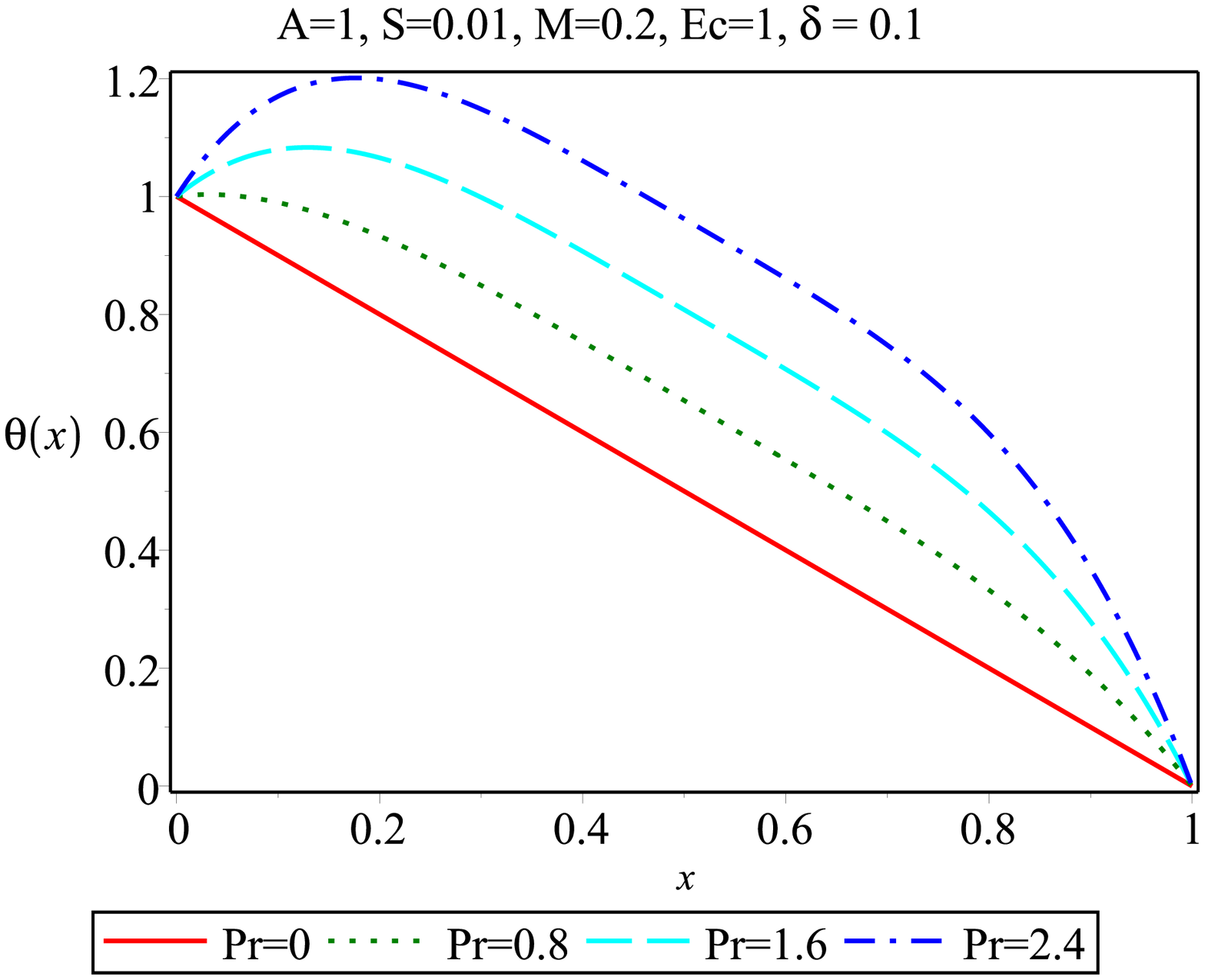}}
\hspace{3mm}
\subfigure[Effect of A on $\theta(x)$]{
\includegraphics*[width=5cm]{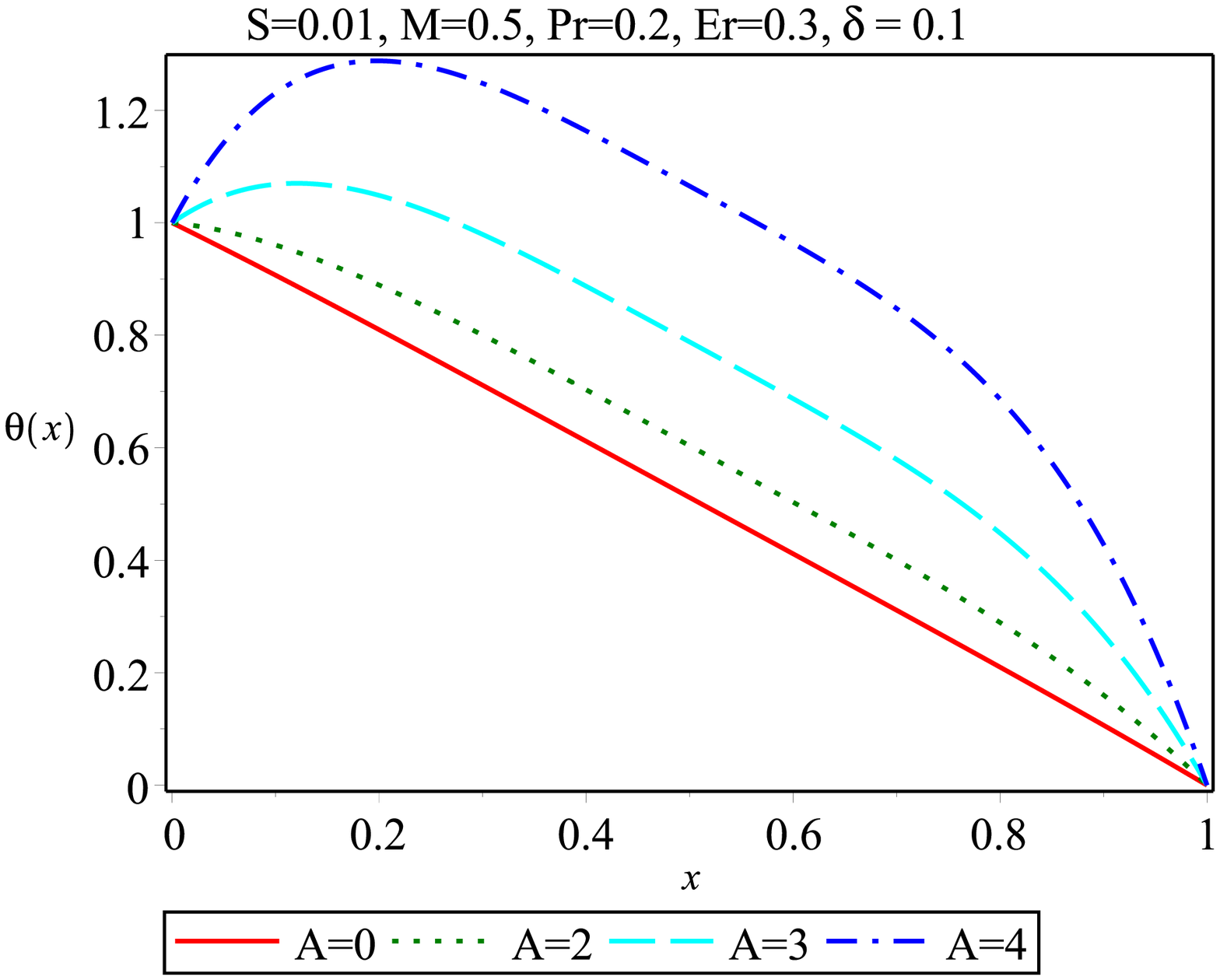}}
\caption{Effect of Pr and A on $\theta(x)$ for the suction cases, Example 1.}
\label{ThetaSuction}
\end{figure}

\begin{figure}[!ht]
\centering
\subfigure[Effect of Pr on $\theta(x)$]{
\includegraphics*[width=5cm]{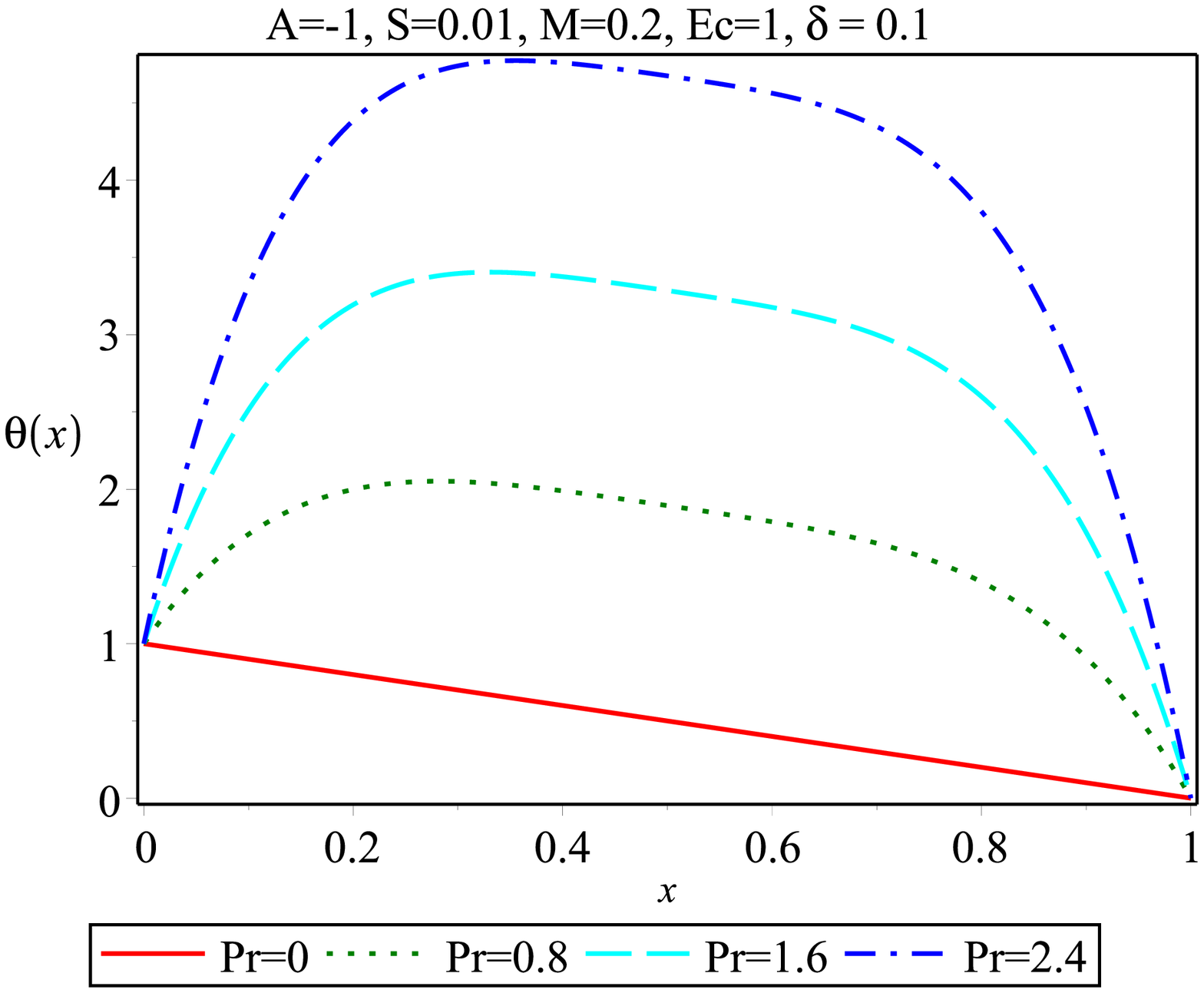}}
\hspace{3mm}
\subfigure[Effect of S on $\theta(x)$]{
\includegraphics*[width=5cm]{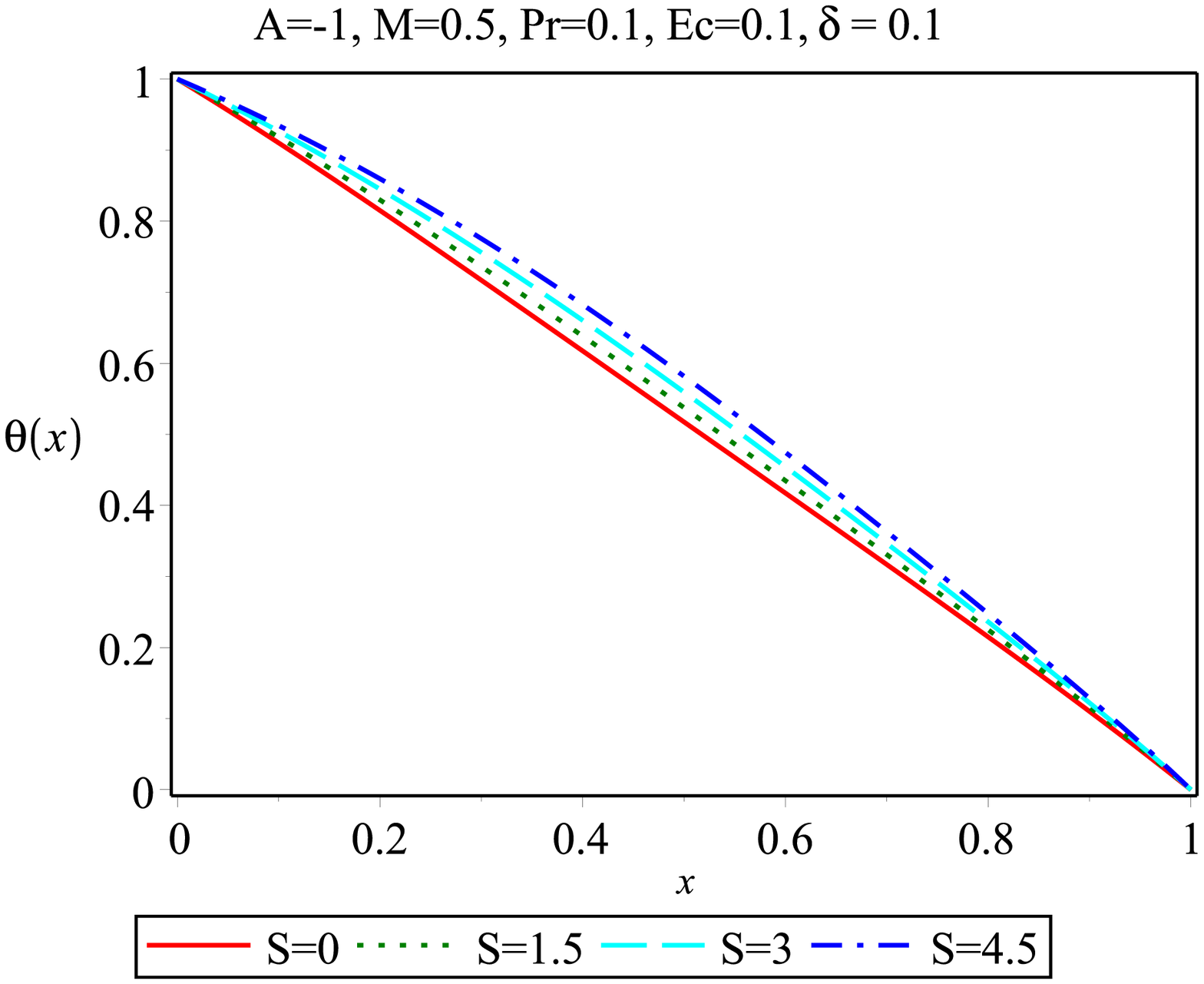}}
\caption{Effect of Pr and S on $\theta(x)$ for the injection cases, Example 1.}
\label{ThetaInjection}
\end{figure}

\begin{figure}[!ht]
\centering
\subfigure[Effect of A on $f'(x)$]{
\includegraphics*[width=5cm]{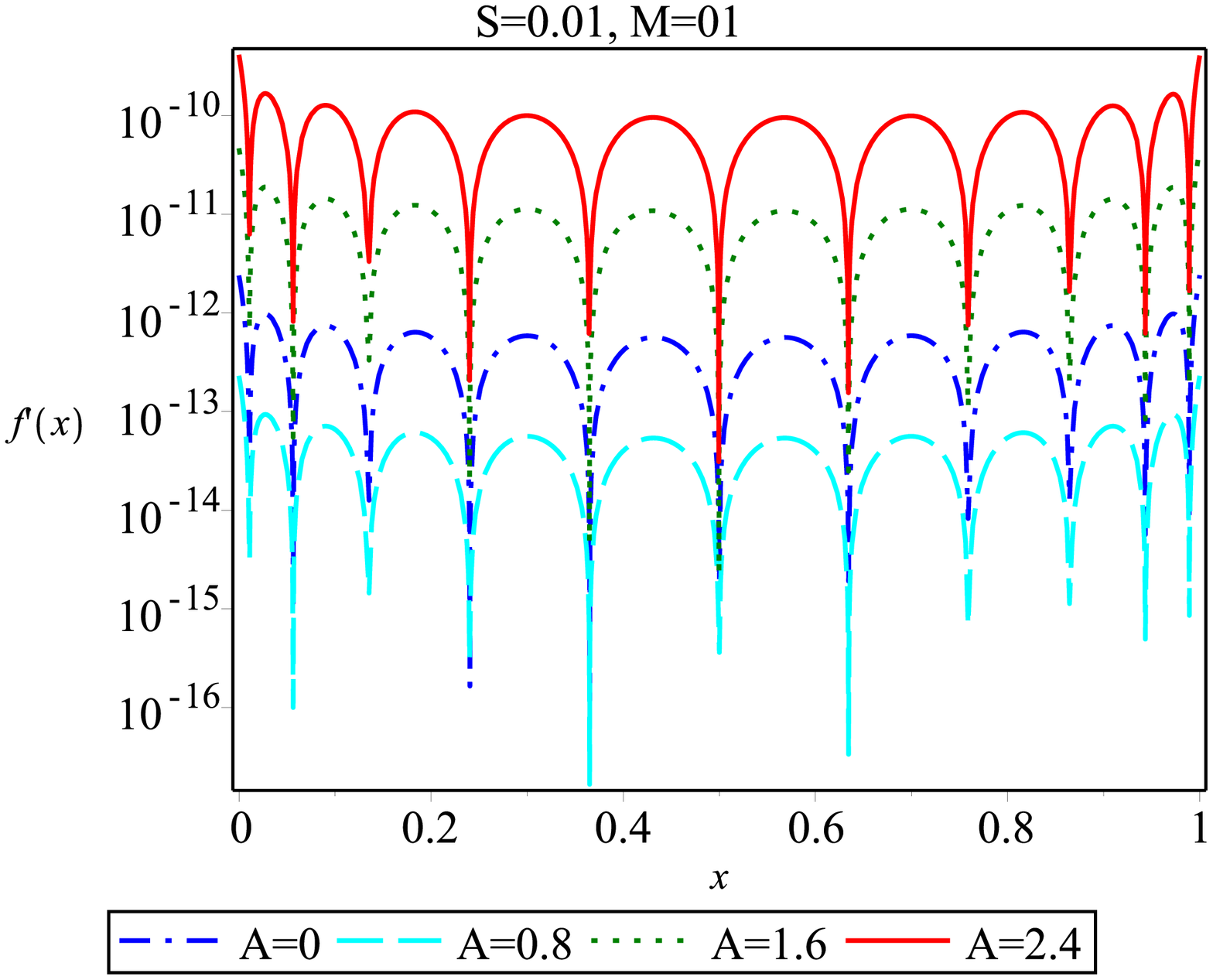}}
\hspace{3mm}
\subfigure[Effect of A on $\theta(x)$]{
\includegraphics*[width=5cm]{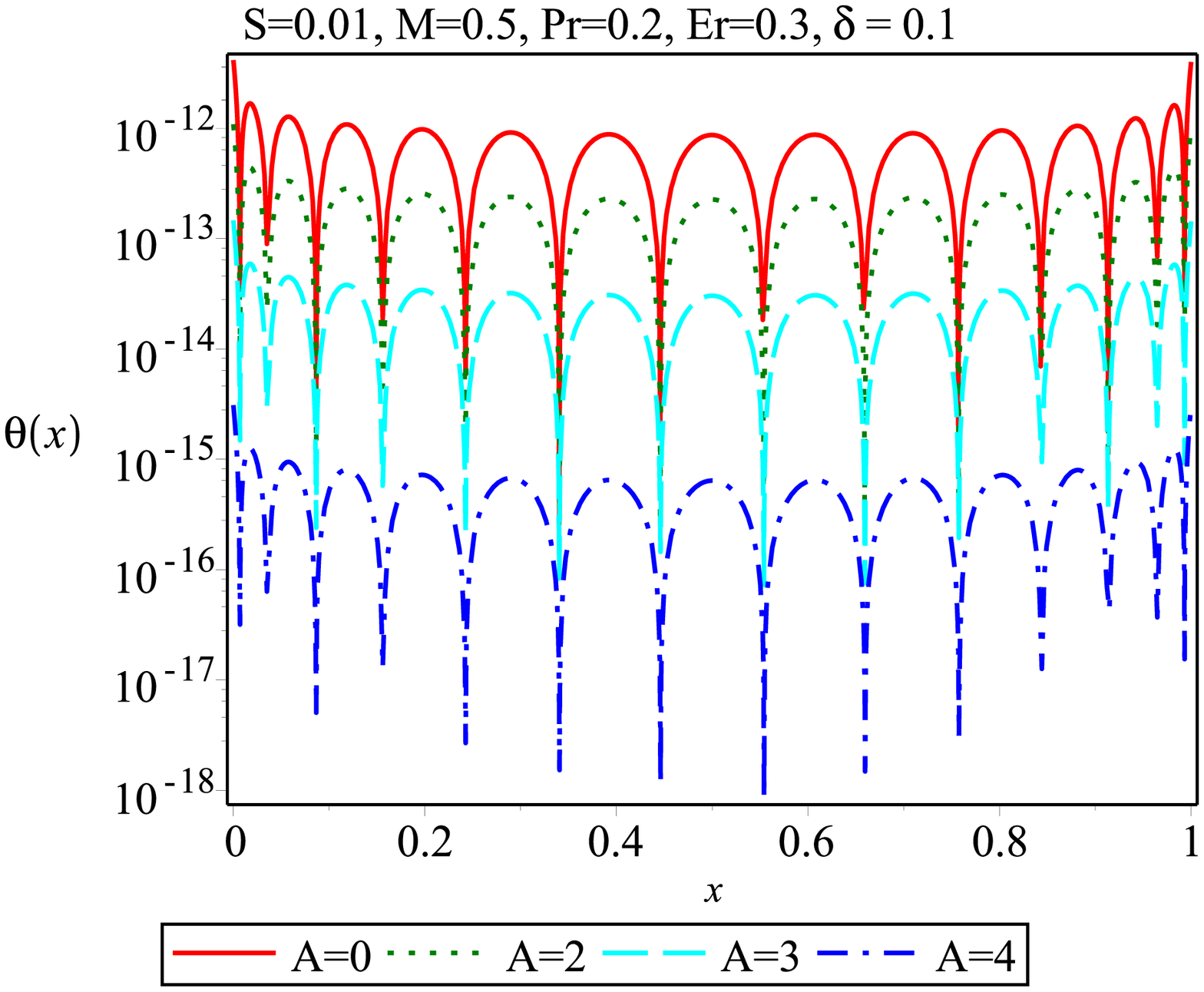}}
\caption{Logarithm graphs of residual error by effects of $A$ on $f'(x)$ and $\theta(x)$ for the suction cases with $N=15$, Example 1.}
\label{ResSuction}
\end{figure}

\begin{figure}[!ht]
\includegraphics[width=7cm]{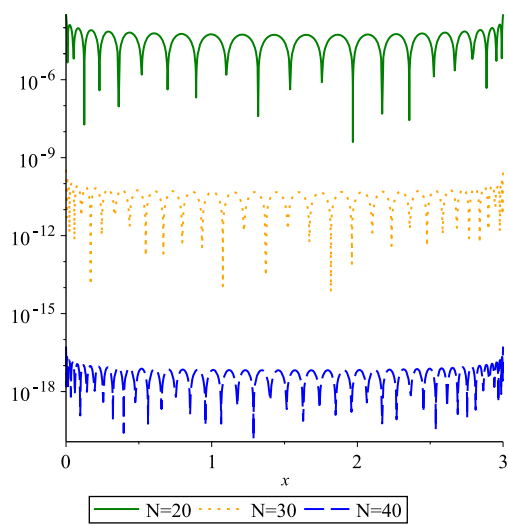}
\centering
\caption{Logarithm graphs of absolute error function, Example 2.}
\label{FigExapmle2}
\end{figure}

\begin{figure}[!ht]
\includegraphics[width=7cm]{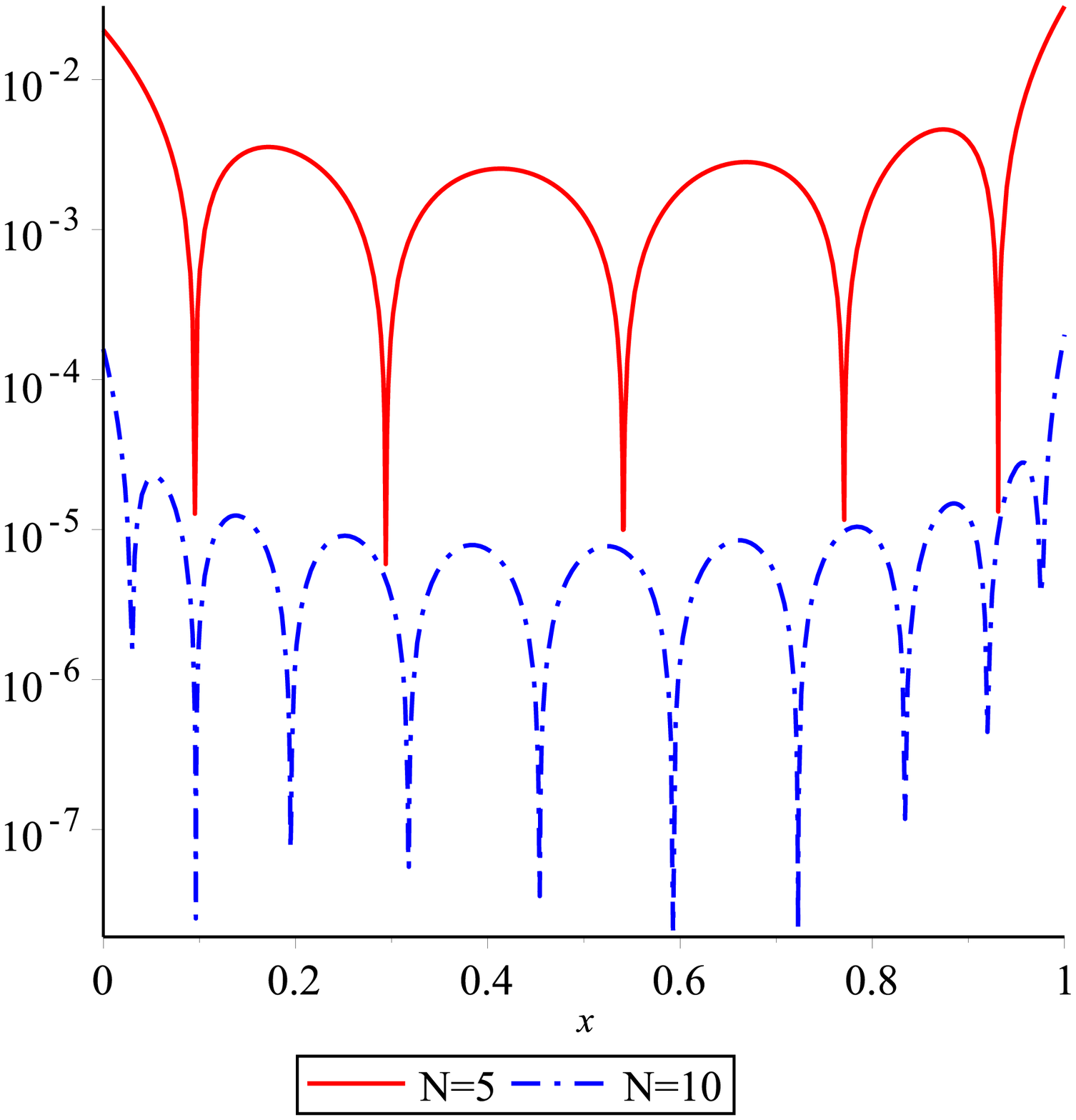}
\centering
\caption{Logarithm graphs of residual error function, Example 3.}
\label{FigExapmle3}
\end{figure}

By using Eqs. (\ref{equation09}) and (\ref{equation16}) we have:
\begin{eqnarray}
\nonumber&&f^{(j)}(x)=F^{T}\textbf{D}^{j}Q(x),\\
\nonumber&&x f^{(j)}= (U^{T}Q(x))(F^{T}\textbf{D}^{j}Q(x))=U^{T}Q(x)Q(x)^{T}\underbrace{(F^{T}\textbf{D}^{j})^{T}}_\text{Z}\simeq U^{T}\tilde{Z}Q(x)\\
\nonumber&&ff''(x)= F^{T}BF^{T}\textbf{D}^{2}J(X)= F^{T}Q(x)Q(x)^{T}(x) \underbrace{(\textbf{D}^{2})^{T}F}_\text{O}\simeq F^{T}\tilde{O}Q(x),\\
\nonumber&&f'^{2}(x)= (F^{T}\textbf{D}Q(x))^{2}=F^{T}\textbf{D}Q(x)Q(x)^{T}\underbrace{D^{T}F}_\text{W}\simeq F^{T}\textbf{D}\tilde{W}Q(x),\\
\nonumber&&f''^{2}(x)= (F^{T}\textbf{D}^{2}Q(x))^{2}=F^{T}\textbf{D}^{2}Q(x)Q(x)^{T}(x)\underbrace{(\textbf{D}^2)^{T}F}_\text{E}\simeq F^{T}\textbf{D}^{2}\tilde{E}Q(x),\\
\nonumber&&\theta^{(j)}(x)= C^{T}\textbf{D}^{j}Q(x),\\
\nonumber&&x\theta'(x)= (U^{T}Q(x))(C^{T}\textbf{D}Q(x))=U^{T}Q(x)Q(x)^{T}(x)\underbrace{\textbf{D}^{T}C}_\text{V}\simeq U^{T}\tilde{V}Q(x),\\
\nonumber&&f\theta'=F^{T}Q(x)C^{T}\textbf{D}Q(x)=F^{T}x^{T}\underbrace{\textbf{D}^{T}C}_\text{H}\simeq F^{T}\tilde{H}Q(x),
\end{eqnarray}
where $Q(x)Q(x)^{T}Z\simeq \tilde{Z}Q(x)$, $Q(x)Q(x)^{T}O\simeq \tilde{O}Q(x)$, $Q(x)Q(x)^{T}E\simeq \tilde{E}Q(x)$, $Q(x)Q(x)^{T}V\simeq \tilde{V}Q(x), Q(x)Q(x)^{T}H\simeq \tilde{H}Q(x)$, and  $\textbf{D}^{j}$ is the $jth$ power of the matrix $\textbf{D}$ and we also have:
\begin{eqnarray}
\nonumber x=X(x)^{T}u=Q(x)^{T}\underbrace{(\textbf{M}^{-1})^{T}u}_\text{U}=Q(x)^{T}U,
\end{eqnarray}
where $u=[0,1,0,0,...,0]^{T}$. Therefore, Eqs. (\ref{equation18}) and (\ref{equation19})  can be rewritten as:
\begin{small}
\begin{eqnarray}
&& F^{T}\textbf{D}^{4}Q(x)-S\Big(U^{T}\tilde{Z}Q(x)+3F^{T}\textbf{D}^{2}Q(x)-2F^{T}\tilde{O}Q(x)\Big)-G^{2}F^{T}\textbf{D}^{2}Q(x)=0,\label{equation23}\\
\nonumber&& C^{T}\textbf{D}^{j}Q(x) +S~Pr\Big(2F^{T}\tilde{H}Q(x)-U^{T}\tilde{V}Q(x)\Big)-Pr~Ec\Big(F^{T}\textbf{D}^{2}\tilde{E}Q(x) +12\delta^2F^{T}\textbf{D}\tilde{W}Q(x) \Big)=0,\\
&&\label{equation24}
\end{eqnarray}
\end{small}
Now, by using Tau method \cite{refExapmle02} and Eq. (\ref{equation04}), we have:
\begin{eqnarray}
\nonumber&&\int_{0}^{1}  \bigg(F^{T}\textbf{D}^{4}-S\Big(U^{T}\tilde{Z}+3F^{T}\textbf{D}^{2}-2F^{T}\tilde{O}\Big)-G^{2}F^{T}\textbf{D}^{2}\bigg)Q(x)Q(x)^{T}{dx}\\
\nonumber&&= \bigg(F^{T}\textbf{D}^{4}-S\Big(U^{T}\tilde{Z}+3F^{T}\textbf{D}^{2}-2F^{T}\tilde{O}\Big)-G^{2}F^{T}\textbf{D}^{2}\bigg)\int_{0}^{1} Q(x)Q(x)^{T}{dx}\\
&&= \bigg(F^{T}\textbf{D}^{4}-S\Big(U^{T}\tilde{Z}+3F^{T}\textbf{D}^{2}-2F^{T}\tilde{O}\Big)-G^{2}F^{T}\textbf{D}^{2}\bigg)K=0,\label{equation25}\\
\nonumber&&\int_{0}^{1}  \bigg(C^{T}\textbf{D}^{j} +S~Pr\Big(2F^{T}\tilde{H}-U^{T}\tilde{V}\Big)-Pr~Ec\Big(F^{T}\textbf{D}^{2}\tilde{E} +12\delta^2F^{T}\textbf{D}\tilde{W} \Big)\bigg)Q(x)Q(x)^{T}{dx}\\
&&=\bigg(C^{T}\textbf{D}^{j} +S~Pr\Big(2F^{T}\tilde{H}-U^{T}\tilde{V}\Big)-Pr~Ec\Big(F^{T}\textbf{D}^{2}\tilde{E} +12\delta^2F^{T}\textbf{D}\tilde{W} \Big)\bigg)K=0.\label{equation26}
\end{eqnarray}
Since $K$ is an invertible matrix, thus one has:
\begin{eqnarray}
&&F^{T}\textbf{D}^{4}-S\Big(U^{T}\tilde{Z}+3F^{T}\textbf{D}^{2}-2F^{T}\tilde{O}\Big)-G^{2}F^{T}\textbf{D}^{2}=0,\label{equation27}\\
&&C^{T}\textbf{D}^{j} +S~Pr\Big(2F^{T}\tilde{H}-U^{T}\tilde{V}\Big)-Pr~Ec\Big(F^{T}\textbf{D}^{2}\tilde{E} +12\delta^2F^{T}\textbf{D}\tilde{W} \Big)=0,\label{equation28}
\end{eqnarray}
 In this method, to satisfy the boundary conditions Eq. (\ref{equation20}), we have substituted the four following equations in four rows of Eq. (\ref{equation27}):
\begin{eqnarray}
&&\nonumber f(0)=F^{T}Q(0)=A,~~~f'(0)=F^{T}\textbf{D}Q(0)=0, \label{equation29} \\
\nonumber&&f(1)=F^{T}Q(1)=\frac{1}{2},~~~f'(1)=F^{T}\textbf{D}Q(1)=0,
\end{eqnarray} 
and we have also substituted two following  equations in two rows of Eq. (\ref{equation28}):
\begin{eqnarray}
&&\nonumber\theta(0)=C^{T}Q(0)=1,\label{equation30} \\
\nonumber&&\theta(1)=C^{T}\textbf{D}Q(1)=0,
\end{eqnarray} 
 Finally, we generate $2N + 2$ algebraic equations, therefore by solving above  equations the unknown vectors $F$ and $C$ are achieved.\\
Effects of deformation parameter S and porosity parameter A are displayed in Figs. (\ref{FSection}) and (\ref{FInjection}).\\
Fig. (\ref{ThetaSuction}) shows the influences of flow parameters (in suction case) on the temperature profile.
Fig. (\ref{ThetaInjection}) shows the influences of physical parameters on velocity and temperature distributions in the case of injection.\\
Fig.(\ref{ResSuction}) shows residual error by the SBT method with $N=15$  on $f'(x)$ and $\theta(x)$ for the suction cases.\\
Tables (1) and  (2) show a comparison between the given values $f'(x)$ and $\theta(x)$ by the variational iteration method (VIM)\cite{refExapmle01} and the SBT methods with $N=15$ for $A=0.1, S=0.1, M=0.2, Pr=0.3, Ec=0.2$ and $\delta=0.1$. It is evident from the table that the present method is efficiency and accuracy of VIM method. Table 3 displays the numerical values of Nusselt number of different values of flow parameters.\\\\
\textbf{Example 2: }
In this example, we consider the following Lane-Emden type equation on the interval $[0,3]$  as follows \cite{refBesselPol13}:
\begin{eqnarray}\label{equation31}
y''(x)+\frac{2}{x}y'(x)-2(2x^2+3)y(x)=0,~~~~0\leqslant x\leqslant3,
\end{eqnarray} 
with the initial conditions
\begin{eqnarray}\label{equation32}
y(0)=1,~~~~~y'(0)=0,
\end{eqnarray} 
The exact solution of this problem is
\begin{eqnarray}\label{equation33}
y(x)=e^{x^{2}}.
\end{eqnarray}
To solve this example, we approximate $y''(x)$  functions by $N+1$ basis of the Bessel polynomials on the interval $[0,3]$ as bellows:
\begin{eqnarray}\label{equation34}
&&y''(x)\simeq \sum_{i=0}^{N} {a_{i}Q_{i}(x)}=A^{T}Q(x)
\end{eqnarray}
where 
$
A=[a_0,a_1,...,a_{N}]^T.
$
By using Eqs. (\ref{equation09}) and (\ref{equation16}) we have:
\begin{eqnarray}
&&\nonumber x=X(x)^{T}w=Q(x)^{T}\underbrace{(\textbf{M}^{-1})^{T}w}_\text{U}=Q(x)^{T}U,\\
&&\nonumber (2x^3+3x)=X(x)^{T}r=Q(x)^{T}\underbrace{(\textbf{M}^{-1})^{T}r}_\text{Z}=Q(x)^{T}Z,\\
&&\nonumber 1=X(x)^{T}p=Q(x)^{T}\underbrace{(\textbf{M}^{-1})^{T}p}_\text{V}=Q(x)^{T}V,
\end{eqnarray}
For approximating the solution $y'(x)$ and $y(x)$,  we apply one and two-time integration on both sides Eq. (\ref{equation34}):

\begin{eqnarray}
&&\nonumber  \int_{0}^{x}y''(x)dx=y'(x)-y'(0)\simeq \int_{0}^{x}A^{T}Q(x)dx \Rightarrow y'(x) \simeq A^{T}\textbf{I}Q(x),\\ \nonumber \\
&&\nonumber y''(x)x=A^{T}Q(x)Q(x)^{T}U\simeq A^{T}\tilde{U}Q(x),\\ \nonumber \\
&&\nonumber  \int_{0}^{x}y'(x)dx=y(x)-y(0)\simeq \int_{0}^{x}A^{T}\textbf{I}Q(x)dx-1\\
&&\nonumber \Rightarrow y(x)=A^{T}\textbf{I}^{2}Q(x)+V^{T}Q(x)=\underbrace{(A^{T}\textbf{I}^{2}+V^{T})}_\text{L}Q(x)= L^{T}Q(x),\\
&&\nonumber y(x)(2x^3+3x)=L^{T}Q(x)Q(x)^{T} Z\simeq L^{T}\tilde{Z}Q(x),
\end{eqnarray}
where $w=[0, 1, 0, 0, ..., 0]^{T}$,  $r=[0, 3, 0, 2, 0, ..., 0]^{T}$ and $p=[1, 0, 0, ..., 0]^{T}$. Therefore, Eq. (\ref{equation31}) can be rewritten as:
\begin{eqnarray} \label{equation35}
&& A^{T}\tilde{U}Q(x)+2A^{T}\textbf{I}Q(x)-2L^{T}\tilde{Z}Q(x)=0,
\end{eqnarray}
As in a typical Tau method \cite{refExapmle02} , we have:
\begin{eqnarray}\label{equation36}
&& A^{T}\tilde{U}+2A^{T}\textbf{I}-2L^{T}\tilde{Z}=0.
\end{eqnarray}
Finally, we generate $N + 1$ algebraic equations, therefore by solving above  equations the unknown vector $A$ is achieved.\\
 Table (4) reports the comparison of $y(x)$ given by exact value, the Hermite functions collocation (HFC) method applied by Parand et al. \cite{Parand.Dehghan.Communications} and  the proposed method  in this paper with $N=40$ . The absolute error graphs of  Lane-Emden equation obtained by present method for $N = 20, 30$ and $40$ are shown in Figure (\ref{FigExapmle2}).\\\\
\textbf{Example 3: }
We consider the following nonlinear Abel equation on the interval $[0,1]$ as follows:
\begin{eqnarray}\label{equation38}
y'(x)=sin(x)y^{3}(x)-xy^{2}(x)+x^{2}y(x)-x^3,~~~~0\leqslant x\leqslant1,
\end{eqnarray} 
with the boundary condition
\begin{eqnarray}\label{equation39}
y(0)=0.
\end{eqnarray}\label{equation40}
To solve this example, we approximate $y'(x)$  functions by $N+1$ basis of shifted Bessel polynomials on the interval $[0,1]$  as bellows:
\begin{eqnarray}\label{equation34}
&&y'(x)\simeq \sum_{i=0}^{N} {a_{i}Q_{i}(x)}=A^{T}Q(x),
\end{eqnarray}
where 
$
A=[a_0,a_1,...,a_{N}]^T.
$
By using Eqs. (\ref{equation09}) and (\ref{equation16}) we have:
\begin{eqnarray}
&&\nonumber x=X(x)^{T}w=Q(x)^{T}\underbrace{(\textbf{M}^{-1})^{T}w}_\text{U}=Q(x)^{T}U,\\
&&\nonumber x^{2}=X(x)^{T}r=Q(x)^{T}\underbrace{(\textbf{M}^{-1})^{T}r}_\text{Z}=Q(x)^{T}Z,\\
&&\nonumber x^{3}=X(x)^{T}p=Q(x)^{T}\underbrace{(\textbf{M}^{-1})^{T}p}_\text{V}=Q(x)^{T}V,\\
&&\nonumber sin(x)=Q(x)^{T}S,\\
&&\nonumber y(x)=A^{T}\textbf{I}Q(x),\\
&&\nonumber y^{2}(x)x=\Big(A^{T}\textbf{I}Q(x)Q(x)^{T}\underbrace{\textbf{I}^{T}A}_\text{W}\Big)Q(x)^{T}U\simeq A^{T}\textbf{I}\tilde{W}\tilde{U}Q(x),\\
&&\nonumber y(x)x^{2}=A^{T}\textbf{I}Q(x)Q(x)^{T}Z\simeq A^{T}\textbf{I}\tilde{Z}Q(x),\\
&&\nonumber y^{3}(x)sin(x)=\Big(A^{T}\textbf{I}\tilde{W}Q(x)Q(x)^{T}\underbrace{\textbf{I}^{T}A}_\text{W}\Big)Q(x)^{T}S\simeq A^{T}\textbf{I}\tilde{W}^{2}\tilde{S}Q(x),
\end{eqnarray} 
where $w=[0, 1, 0, 0, ..., 0]^{T}$,  $r=[0, 0, 1, 0, ..., 0]^{T}$, $p=[0, 0, 0, 1, 0,  ..., 0]^{T}$ and
$S^{T}=\langle sin(x),Q^{T}(x)\rangle_{w}\langle Q(x),Q^{T}(x)\rangle_{w}^{-1},$.\\
Therefore, Eq. (\ref{equation38}) can be rewritten as:
\begin{eqnarray}
A^{T}Q(x) - A^{T}\textbf{I}\tilde{W}^{2}\tilde{S}Q(x) + A^{T}\textbf{I}\tilde{Z}Q(x) - A^{T}\textbf{I}\tilde{Z}Q(x) + V^{T}Q(x)=0,
\end{eqnarray} 
As in a typical Tau method \cite{refExapmle02} , we have:
\begin{eqnarray}
A^{T} - A^{T}\textbf{I}\tilde{W}^{2}\tilde{S} + A^{T}\textbf{I}\tilde{Z} - A^{T}\textbf{I}\tilde{Z} + V^{T}=0,
\end{eqnarray} 
In this step, we generate $N + 1$ algebraic equations, therefore by solving above  equations the unknown vector $A$ is achieved.\\
 The residual error graphs of the nonlinear Abel equation obtained by present method for $N = 5$ and $10$ are shown in Figure (\ref{FigExapmle3}). Table (5) displays given values $y(x)$  by the fractional order of the Chebyshev orthogonal functions (FCFs) collocation method used by parand et al. \cite{Delkhosh.Nikarya}, the present method and the residual error with $N=10$.\\\\
\textbf{Example 4: }
Consider the nonlinear the standard Lane-Emden problem
\begin{eqnarray}\label{standardLaneEmden}
y''(x)+\frac{1}{x}y'(x)+y^{2}=0,~~~~0\leqslant x\leqslant2
\end{eqnarray}
with the boundary condition
\begin{eqnarray}
y(0)=1,~~~y'(0)=0.
\end{eqnarray}

Thus, the solution of the Lane-Emden equation is obtained by the SBT method with  $N=8$ and $N=12$.

\begin{figure}[!ht]
\includegraphics[width=7cm]{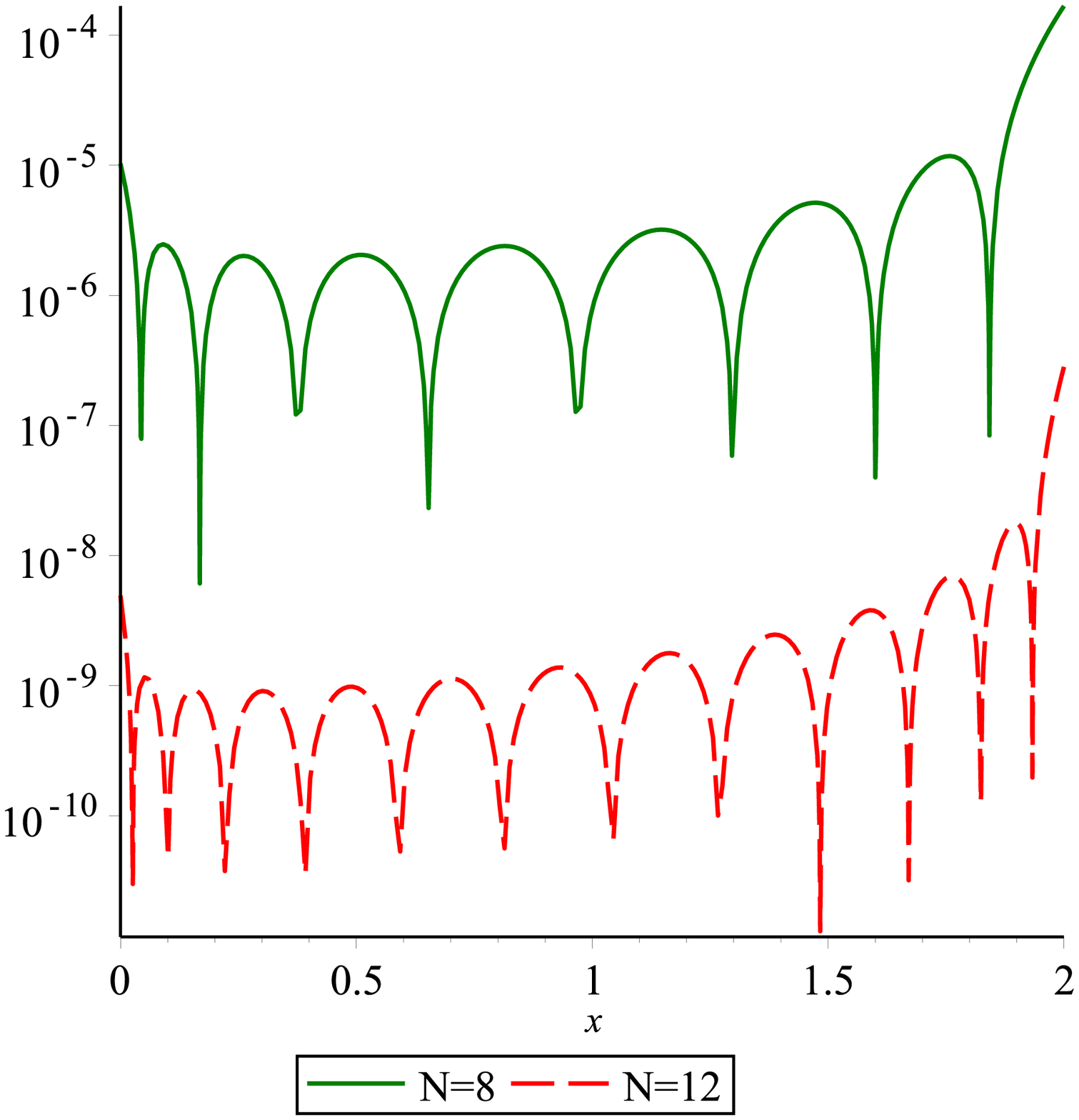}
\centering
\caption{Logarithm graphs of residual error function, Example 4.}
\label{FigExapmle4}
\end{figure}

 The residual error graphs of the nonlinear standard Lane-Emden  obtained by the present method for $N = 8$ and $12$ are shown in Figure (\ref{FigExapmle4}). Table (6) reports  the comparison of the obtained values $y(x)$  by the present method with $N=12$, Horedt \cite{refExapmle04}, the Bessel orthogonal functions collocation (BFC) method utilized by Parand et al. \cite{refBesselPol05} and the Chebyshev orthogonal functions (GFCFs) collocation method proposed by Parand et al. \cite{Parand.Delkhosh.Teknoligy}.\\\\
\textbf{Example 5: }
ّFinally, consider the boundary-value the Troesch's problem
\begin{eqnarray}
y''(x)=\gamma \sinh(\gamma y(x)), ~~~~0 \leqslant x \leqslant 1,
\end{eqnarray}
with the boundary conditions
\begin{eqnarray}\label{boundary.example5}
y(0)=0,~~~y(1)=1,
\end{eqnarray}
We approximate $y(x)$  functions by $N+1$ basis of shifted Bessel polynomials on the interval $[0,1]$  as bellows:
\begin{eqnarray}\label{equation34}
&&\nonumber y(x)\simeq \sum_{i=0}^{N} {a_{i}Q_{i}(x)}=A^{T}Q(x)
\end{eqnarray}
where $A=[a_0,a_1,...,a_{N}]^T$ . 

By using the SBT technique described in the paper, we can generate a set of linear algebraic equations as follows:
\begin{eqnarray}\label{algebraic.example5}
A^{T}\textbf{D}^{2}-\gamma(\gamma A^{T}+\frac{\gamma^{3}}{3!}A^{T}\tilde{A}^{2}+\frac{\gamma^{5}}{5!}A^{T}\tilde{A}^{4})=0,
\end{eqnarray}
where $Q(x)Q(x)^{T}A\simeq \tilde{A}Q(x)$.

To satisfy the boundary conditions Eq. (\ref{boundary.example5}), we have substituted the following equations in two rows of Eq. (\ref{algebraic.example5}):
\begin{eqnarray}
&& y(0)=A^{T}Q(0)=0,~~~y(1)=A^{T}Q(0)=1.
\end{eqnarray} 

Table (7) displays  the comparison of the obtained values $y(x)$  by the modified non-linear shooting method (MNLSM) applied by Alias \cite{Alias.Teknology}, the numerical scheme based on the modified homotopy perturbation technique used by \cite{Feng.Computation} and the present method with $N=10$ and $\gamma=0.5$.

\section{Conclusion}
The shifted Bessel Tau (SBT) method is a useful technique for solving ordinary differential equations. Numerical programs using this method are often considerably faster with greater accuracy than other standard techniques such as the variational iteration method, shooting method and homotopy perturbation technique. Interest in the Tau method, for a long time regarded only as a tool for the construction of accurate approximations of a very restricted class of functions, has been enhanced. In this investigation, the Bessel polynomials operational matrices of integration, differentiation, and product are derived. A general procedure for the formulation of these matrices is given. By utilizing the operational matrices  and  the Tau method, we reduce the solution of the original equations to the solution of a nonlinear algebraic equation system. Therefore, all of these equations can be solved by Newton method for the unknown coefficients. Some problems with initial or boundary conditions  are considered in one dimensional contexts, although the proposed technique can be easily implemented for two and three dimensional problems. The main advantage of the present method is its simplicity and convenience for computer algorithms. Illustrative examples demonstrate the validity and applicability of the present method.
%\begin{acknowledgements}
%\end{acknowledgements}


%\footnotesize

%\normalsize

%%%%%%%%%%%%%%%%%%%%%%%%%%%%%%%%%%

%%%%%%%%%%%%%%%%%%%%%%%%%%%%%%%%%%%%%%%%%%

%%%%%%%%%%%%%%%%%%%%%%%%%%%%%%%%%%%%%%%%%%

%%%%%%%%%%%%%%%%%%%%%%%%%%%%%%%%%%%%%%%%%%

%%%%%%%%%%%%%%%%%%%%%%%%%%%%%%%%%%%%%%%%%%

%%%%%%%%%%%%%%%%%%%%%%%%%%%%%%%%%%%%%%%%%%

%%%%%%%%%%%%%%%%%%%%%%%%%%%%%%%%%%%%%%%%%%

%%%%%%%%%%%%%%%%%%%%%%%%%%%%%%%%%%%%%%%%%%

%%%%%%%%%%%%%%%%%%%%%%%%%%%%%%%%%%%%%%%%%%

%%%%%%%%%%%%%%%%%%%%%%%%%%%%%%%%%%%%%%%%%%

%%%%%%%%%%%%%%%%%%%%%%%%%%%%%%%%%%%%%%%%%%

%%%%%%%%%%%%%%%%%%%%%%%%%%%%%%%%%%%%%%%%%%
\footnotesize 
\begin{center}
\begin{table}[ht!]
\label{TblF}
\caption{Comparison of the obtained values $f'(x)$ by the variational iteration method (VIM)\cite{refExapmle01} and the present method, Example 1.}
\centering
%\footnotesize 
\begin{tabular}[ht!]{|c|c|c|c|}\hline
   $x$ &Bessel Tau method&VIM\cite{refExapmle01}&$Res(x)$   \\ \hline

0.2&0.384801280290557 &0.384801&1.39535e-12\\
0.4&0.575554143054330 &0.575554&8.21911e-14\\
0.6&0.575174675564395 &0.575174&1.90287e-12\\
0.8&0.384040432902588 &0.384040&1.91562e-12\\\hline
	
\end{tabular}
\end{table}
\end{center}
%%%%%%%%%%%%%%%%%%%%%%%%%%%%%%%%%%%%%%%%%%

%%%%%%%%%%%%%%%%%%%%%%%%%%%%%%%%%%%%%%%%%%
\footnotesize 
\begin{center}
\begin{table}[ht!]
\label{TblTheta}
\caption{Comparison of the obtained values $\theta(x)$ by the variational iteration method (VIM)\cite{refExapmle01} and the present method, Example 1.}
\centering
%\footnotesize 
\begin{tabular}[ht!]{|c|c|c|c|}\hline

   $x$ &Bessel Tau method&VIM\cite{refExapmle01}&$Res(x)$   \\ \hline

0.2&0.806144282850332 &0.806144&8.25970e-14\\
0.4&0.607022034290869 &0.607022&7.26187e-14\\
0.6&0.407003862839112 &0.407004&7.24240e-14\\
0.8&0.206120555470330 &0.206121&8.19659e-14\\\hline
	
\end{tabular}
\end{table}
\end{center}
%%%%%%%%%%%%%%%%%%%%%%%%%%%%%%%%%%%%%%%%%%

%%%%%%%%%%%%%%%%%%%%%%%%%%%%%%%%%%%%%%%%%%
\footnotesize 
\begin{center}
\begin{table}[ht!]
\label{TblTheta}
\caption{Comparison of values of Nusselt number $ (1-at)^{1/2}Nu$ for different values of $Pr, Ec$, and $\delta$ when $A=0.1,~S=0.1,~M=0.2$ by the present method with $N=15$ and the variational iteration method (VIM) \cite{refExapmle01}, Example 1.}
\centering
%\footnotesize 
\begin{tabular}{|c|c|c|c|c|}\hline
Pr & Ec & $\delta$ & $(1-at)^{1/2}Nu$ & VIM \cite{refExapmle01}\\\hline
0.0 & 0.2 & 0.1& 1.00000000000000000 & 1.00000\\
0.1 &      &        & 1.01893685003010788 & 1.01894 \\
0.2 &      &        & 1.03787156909761297 & 1.03787  \\
0.3 & 0.0 &      &  0.99859957004178043 & 0.99860 \\
      & 0.2 &      &  1.05680415677248143& 1.05680 \\
      & 0.4 &      &  1.11500874350318244 & 1.11501 \\
      & 0.3 & 0.0 & 1.08487154864359339 & 1.08487  \\
      &      & 0.5 &  1.11074408599955706 & 1.11074 \\
      &      & 1.0 &  1.18836169806744806 & 1.18836 \\\hline
\end{tabular}
\end{table}
\end{center}
%%%%%%%%%%%%%%%%%%%%%%%%%%%%%%%%%%%%%%%%%%

%%%%%%%%%%%%%%%%%%%%%%%%%%%%%%%%%%%%%%%%%%

%\begin{figure}[!ht]
%\includegraphics[width=10cm]{Example2y.eps}
%\centering
%\caption{Lane-Emden graph obtained by the present method with $N=20$, Example 2.}
%\label{yExapmle2}
%\end{figure}
%%%%%%%%%%%%%%%%%%%%%%%%%%%%%%%%%%%%%%%%%%

%%%%%%%%%%%%%%%%%%%%%%%%%%%%%%%%%%%%%%%%%%

%%%%%%%%%%%%%%%%%%%%%%%%%%%%%%%%%%%%%%%%%%

%%%%%%%%%%%%%%%%%%%%%%%%%%%%%%%%%%%%%%%%%%
\footnotesize 
\begin{center}
\begin{table}[ht!]
\label{TblF}
\caption{Comparison of the obtained values $y(x)$  by  exact value, HFC method \cite{Parand.Dehghan.Communications} and the present method with $N=40$, Example 2.}
\centering
%\footnotesize 
\begin{tabular}[ht!]{|c|c|c|c|c|}\hline
   $x$& Exact value&HFC \cite{Parand.Dehghan.Communications}&Peresent method  \\\hline

0.01&1.00010000500016667083&1.0000999826&1.00010000500016665722\\
0.02&1.00040008001066773341&1.0004000642&1.00040008001066774881\\
0.05&1.00250312760579508497&1.0025031064&1.00250312760579507309\\
0.10&1.01005016708416805754&1.0100501492&1.01005016708416805546\\
0.20&1.04081077419238822675&1.0408107527&1.04081077419238822399\\
0.50&1.28402541668774148407&1.2840253862&1.28402541668774147818\\
0.70&1.63231621995537897012&1.6323161777&1.63231621995537896294\\
0.80&1.89648087930495135334&1.8964808279&1.89648087930495136037\\
0.90&2.24790798667647141917&2.2479078937&2.24790798667647141232\\
1.00&2.71828182845904523536&2.7182819166&2.71828182845904524184\\
1.5&9.48773583635852572055 &           ---       &9.48773583635852572300\\
2.0&54.5981500331442390781 &           ---       &54.5981500331442390719\\
2.5&518.012824668342025939 &           ---       &518.012824668342025947\\
3.0&8103.08392757538400770 &           ---       &8103.08392757538400765\\ \hline
	
\end{tabular}
\end{table}
\end{center}
%%%%%%%%%%%%%%%%%%%%%%%%%%%%%%%%%%%%%%%%%%

%%%%%%%%%%%%%%%%%%%%%%%%%%%%%%%%%%%%%%%%%%

%%%%%%%%%%%%%%%%%%%%%%%%%%%%%%%%%%%%%%%%%%

%%%%%%%%%%%%%%%%%%%%%%%%%%%%%%%%%%%%%%%%%%
\footnotesize 
\begin{center}
\begin{table}[ht!]
\label{TblF}
\caption{Numerical given values $y(x)$ by FCFs collocation method \cite{Delkhosh.Nikarya} and the present method with $N=10$, Example 3.}
\centering
%\footnotesize 
\begin{tabular}[ht!]{|c|c|c|c|c|}\hline
   $x$ &FCFs \cite{Delkhosh.Nikarya}&Present method&Residual error  \\ \hline

0.1&-2.500e-5&-2.541201381e-5 &2.13471e-6\\
0.2&-4.004e-4&-4.000957821e-4 &1.42107e-6\\
0.3&-2.032e-3&-2.033184449e-3 &3.48798e-6\\
0.4&-6.459e-3&-6.459447204e-3 &7.35610e-6\\
0.5&-1.591e-2&-1.591376982e-2 &6.64911e-6\\
0.6&-3.346e-2&-3.346334651e-2 &1.36557e-6\\
0.7&-6.327e-2&6.327405918e-2 &4.86862e-6\\
0.8&-1.111e-1&-1.111347323e-1 &9.41866e-6\\
0.1&-1.855e-1&1.855105463e-1 &1.23827e-5\\
1.0&-2.999e-1&-2.999554921e-1 &1.98917e-4\\\hline
	
\end{tabular}
\end{table}
\end{center}
%%%%%%%%%%%%%%%%%%%%%%%%%%%%%%%%%%%%%%%%%%

%%%%%%%%%%%%%%%%%%%%%%%%%%%%%%%%%%%%%%%%%%

%%%%%%%%%%%%%%%%%%%%%%%%%%%%%%%%%%%%%%%%%%

%%%%%%%%%%%%%%%%%%%%%%%%%%%%%%%%%%%%%%%%%%
\footnotesize 
\begin{center}
\begin{table}[ht!]
\label{TblF}
\caption{Comparison of the obtained values $y(x)$  by the present method with $N=12$, Horedt \cite{refExapmle04}, BFC method \cite{refBesselPol05} and GFCFs collocation method \cite{Parand.Delkhosh.Teknoligy}, Example 4.}
\centering
%\footnotesize 
\begin{tabular}[ht!]{|c|c|c|c|c|c|}\hline
   $x$& Present method &Horedt \cite{refExapmle04} &BFC \cite{refBesselPol05}& GFCFs \cite{Parand.Delkhosh.Teknoligy}\\\hline

0.1&0.998334998549872 &0.9983350&0.99833499854&0.99833499986\\
0.3&0.985133946938390 &--&--&--\\
0.5&0.959352715810926 &0.9593527&0.95935271580&0.95935271585\\
0.7&0.922170348514590 &--&--&--\\
1.0&0.848654111411546 &0.8486541&0.84865411140&0.84865409603\\
1.5&0.695367147241325 &--&--&--\\
2.0&0.529836429310169 &--&--&--\\\hline
	
\end{tabular}
\end{table}
\end{center}
%%%%%%%%%%%%%%%%%%%%%%%%%%%%%%%%%%%%%%%%%%

%%%%%%%%%%%%%%%%%%%%%%%%%%%%%%%%%%%%%%%%%%
\footnotesize 
\begin{center}
\begin{table}[ht!]
\label{TblF}
\caption{Comparison of the obtained values $y(x)$  by Alias \cite{Alias.Teknology}, Feng \cite{Feng.Computation} and the present method with $N=10$ and $\gamma = 0.5$, Example 5.}
\centering
%\footnotesize 
\begin{tabular}[ht!]{|c|c|c|c|c|c|}\hline
   $x$ &Alias \cite{Alias.Teknology} &Feng \cite{Feng.Computation}&  Present method & Residual eroor\\\hline

0.1&0.09597247 &0.0959477541&0.095944350620621&3.18e-13\\
0.2&0.19218506 &0.1921352537&0.192128750320282&6.98e-12\\
0.3&0.28887905 &0.2888034214&0.288794404891654&1.04e-10\\
0.4&0.38629807 &0.3861955524&0.386184851707410&7.44e-10\\
0.5&0.48441684 &0.4845585473&0.484547171441282&3.63e-09\\
0.6&0.58428140 &0.5841442013&0.584133256467667&1.49e-08\\
0.7&0.68525684 &0.6852105701&0.685201157498259&5.56e-08\\
0.8&0.78807945 &0.7880234321&0.788016532411673&1.91e-07\\
0.9&0.89292601 &0.8928578710&0.892854224345211&6.07e-07\\\hline
	
\end{tabular}
\end{table}
\end{center}
%%%%%%%%%%%%%%%%%%%%%%%%%%%%%%%%%%%%%%%%%%

\end{document}